%% file: 1999-001.tex
\documentclass[11pt]{article}
\usepackage{amsmath}
\usepackage{amsthm}
\usepackage{psfig}
\usepackage{fullpage}

\newcommand {\norm} [1] { \lVert #1 \rVert}
\newcommand {\Norm} [1] {\Bigl\lVert #1 \Bigr\rVert}

\newcommand {\abs} [1] {\lvert #1 \rvert}
\newcommand {\Abs} [1] {\bigl\lvert #1 \bigr\rvert}

\newcommand {\inv} [1] {#1^{-1}}

\newcommand {\ie} {i.e.}

\newcommand{\floor}[1]{\lfloor #1\rfloor}
\newtheorem {thm} {Theorem}[section]
\newtheorem {prop}[thm] {Proposition}
\newtheorem {lem} [thm] {Lemma}
\newtheorem {cor} [thm]  {Corollary}

\theoremstyle {definition}
\newtheorem{defn} {Definition} 
\newtheorem{exmp} {Example} [section]

\theoremstyle {remark}

\numberwithin{equation}{section}

\title{A Conditioning Function for the
Convergence of Numerical ODE Solvers and Lyapunov's Theory of
Stability}

\author {Divakar Viswanath \thanks{ 
Research at MSRI is supported
in part by NSF grant DMS-9701755}}
\date{22 December 1998}
\begin{document}
\maketitle
\begin{abstract}
For the ordinary differential equation (ODE) $\dot{x}(t) = f(t,x)$,
$x(0) = x_0$, $t\geq 0$, $x\in R^d$, assume $f$ to be at least
continuous in $t$ and locally Lipshitz in $x$, and if necessary,
several times continuously differentiable in $t$ and $x$.
We associate a conditioning function $E(t)$ with each solution $x(t)$
which captures the accumulation of global error in a numerical
approximation in the following sense: if $\tilde{x}(t;h)$ is an 
approximation derived from a single step method of time step $h$ and
order $r$ then $\norm{\tilde{x}(t;h) - x(t)} < K(E(t)+\epsilon)h^r$
for $0\leq t\leq T$, any $\epsilon > 0$, sufficiently small $h$,
and a constant $K>0$.

Using techniques from the stability theory of differential equations,
this paper gives conditions on $x(t)$ for $E(t)$ to be upper bounded
linearly or by a constant for $t\geq 0$. More concretely, these
techniques give constant or linear bounds on $E(t)$ when $x(t)$ is
a trajectory of a dynamical system which falls into a stable, hyperbolic
fixed point; or into a stable, hyperbolic cycle; or into a normally hyperbolic 
and contracting manifold with quasiperiodic flow on the manifold.
\end{abstract}

\input s0.ltx

\input s1.ltx

\input s2.ltx

\input s3.ltx

\input s4.ltx
\input s5.ltx
\input s6.ltx
\input s7.ltx

\input bib.ltx
\noindent
Mathematical Sciences Research Institute (MSRI)\\
1000 Centennial Drive\\
Berkeley, CA 94720.\\
divakar@cs.cornell.edu

\end{document}

%% file: s0.ltx
\section{Introduction}

For the system of ordinary differential equations $\dot{x}(t) =
f(t,x)$, $t\geq 0$, $x\in R^d$, the initial value problem $x(t_0) =
x_0$, $t_0\geq 0$, may not have a solution, and when the solution
exists it may not be unique. However, if $f(t,x)$ is continuous in $t$
and not only continuous but also locally Lipshitz in $x$, there is a
unique solution of the initial value problem for $t_0 \leq t < t_0 +
\epsilon$, $\epsilon > 0$, which we  denote by $x(t; t_0,
x_0)$. We make these assumptions about $f(t,x)$ and  only consider
solutions $x(t;t_0,x_0)$ that can be continued till $t=\infty$.
Usually $t_0=0$, and we denote these solutions by $x(t;x_0)$ or simply
$x(t)$.

Even at this level of generality, one might ask how accurately
$x(t;x_0)$ can be approximated by a numerical method. A reasonable first
guess is to look at the $d\times d$ matrix
\begin{equation*}
E'(t) = \frac{\partial x(t;x_0)}{\partial x_0},
\end{equation*}
for $0\leq t<\infty$; since $E'(t)$ gives the sensitivity of $x(t;x_0)$
with respect to $x_0$, one might hope to relate it to the accumulation
of global error. However, numerical methods introduce discretization error
not just at $t=0$ but at every time step of the integration, and
therefore, $E'(t)$ proves to be an insufficient concept.

The following conditioning function $E(t)$ associated with $x(t;x_0)$,
or more briefly $x(t)$, does capture the accumulation of global error:
\begin{equation}
E(t) = \sup_{v(s)}\Norm{\int_0^t \frac{\partial x(t)}{\partial x(s)} 
       v(s) ds},
\label{eqn-0-1} 
\end{equation}
with the supremum taken over continuous functions $v:[0,t] \rightarrow
R^d$ with $\norm{v(s)} \leq 1$ for $0 \leq s \leq t$. Unlike $E'(t)$,
$E(t)$ takes into account the sensitivity of $x(t)$ with respect to
all $x(s)$, $0 \leq s \leq t$. In the definitions of $E(t)$ and
$E'(t)$ above, we have assumed $f(t,x)$ to be continuously
differentiable with respect to $t$ and $x$.  All vector norms in this
paper are Euclidean norms and all matrix norms are the corresponding
induced norms.

Let $\tilde{x}(t;x_0;h)$ denote the approximation to $x(t;x_0)$
computed by a single step method of order $r$. The numerical method
gives $\tilde{x}(t;x_0;h)$ at all positive integer multiples of $h$.
For $kh < \tau < (k+1)h$, $k=0,1,2,\ldots$, we define
$\tilde{x}(\tau;x_0;h)$ by following the solution exactly from the
initial point $t=kh$ and $x=\tilde{x}(kh;x_0;h)$ till
$t=\tau$. Therefore, $\tilde{x}(t;x_0;h)$ can be discontinuous only at
$t=kh$, $k=1,2,\ldots$ The magnitude of this discontinuity
\begin{equation*} 
\norm{\tilde{x}((k+1)h;x_0;h) - x((k+1)h; kh, \tilde{x}(kh;x_0;h))}
\end{equation*}
is the local discretization error at $t=(k+1)h$. We write this
discretization error as $K_{k+1} h^{r+1} v_{k+1}$, where $v_{k+1}\in
R^d$ with $\norm{v_{k+1}} = 1$ and $K_{k+1}$ is a non-negative real
number. Thus the direction of the local discretization error at the
$i$th step is $v_i$ and its magnitude is $K_i h^{r+1}$. 

Let us assume that the magnitude of the local discretization error
is bounded above by $Kh^{r+1}$ for a constant $K>0$. This assumption
can hold for example for a Runge-Kutta method of order $r$ if
$f(t,x)$ is $r+1$ times continuously differentiable with respect to
$t$ and $x$; we discuss this assumption further in Section 2. Then
by Theorem \ref{thm-3-3}, given $\epsilon > 0$ and $T>0$,
\begin{equation}
\norm{x(t;x_0) - \tilde{x}(t;x_0;h)} < (E(t)+\epsilon)Kh^r
\label{eqn-0-2}
\end{equation}
for $0\leq t\leq T$ and sufficiently small $h$. In this bound on the
global error, $E(t)$, which is given by \eqref{eqn-0-1}, is
independent of details of the numerical method. Further, a bound like
\eqref{eqn-0-2} will not hold for any real-valued function of $t$
which is strictly less than $E(t)$ at some value of $t$. It is
because of these reasons that we call $E(t)$ a {\it conditioning
function}.

Let us now draw an analogy between $E(t)$ and the absolute condition 
number of a multivariate function $g(x)$. The formula (see \cite{TB})
\begin{equation*}
\lim_{\delta \rightarrow 0^+} \sup_{\delta x < \delta}
\frac{\norm{g(x+\delta x) - g(x)}}{\delta}
\end{equation*}
is the analogue of \eqref{eqn-0-1}; this absolute condition number 
measures the sensitivity of $g(x)$ with respect to small changes in
$x$. For a stable numerical method, the error in evaluating $g(x)$ 
is governed by this absolute condition number in a manner similar
to the dependence of global error on $E(t)$ given by 
\eqref{eqn-0-2}.

Sections 2,3, and 4 define $E(t)$ and develop its properties. We do
not begin with \eqref{eqn-0-1} as the definition of $E(t)$. Section 2
introduces a model of discretization errors which is similar to but
more general than the model of Stuart and Humphries \cite{SH}. Section
3 defines $E(t)$ using this model. The expression for $E(t)$ in
\eqref{eqn-0-1} is derived in Corollary \ref{cor-5-0} of Section
6. Let us mention the formal similarity of \eqref{eqn-0-1} to bounds
on the global error derived by Iserles and S\"{o}derlind \cite{IS}
using the Alexseev-Gr\"{o}bner lemma. The results in this half of the
paper argue that $E(t)$ is the appropriate vehicle for a study of
global errors.

From Section 5 onwards, we undertake the task of relating $E(t)$ to
stability properties of $x(t;x_0)$. The relation of $E(t)$ to the
accumulation of global error as $t$ increases is clear enough from
\eqref{eqn-0-2}. If $E(t)$ is bounded by a constant or linearly or by
a polynomial of low degree in $t$, the accurate approximation of the
trajectory $x(t;x_0)$ can be considered a tractable problem; but if
$E(t)$ increases exponentially in $t$ accurate approximation of
$x(t;x_0)$ is pretty intractable. The standard technique
of bounding global errors using the Lipshitz constant is used in
convergence proofs of numerical methods \cite{HNW}. But the bounds on
$E(t)$ obtained this way increase exponentially in $t$ and are of
hardly any other use.

We upper bound $E(t)$ by making stability assumptions on $x(t;x_0)$. 
Sections 5 and 6 derive linear and constant bounds on $E(t)$ by making
stability assumptions on $x(t;x_0)$. This connection to stability is
somewhat subtle; there are exponentially stable examples with
exponentially increasing $E(t)$. However, the work of researchers
who followed Lyapunov allows us to clarify and circumvent the difficulties
in relating $E(t)$ to stability theory. See Table \ref{table-1} 
in Section 6 for 
a summary.

Let us now make some prefatory remarks about the stability theory
of ordinary differential equations. The theory of stability of ordinary differential equations was initiated
by A.M. Lyapunov in 1892 \cite{Lyapunov}. One stream of research which
emanates from this remarkable work  is about first approximations.
Let $x(t)=x(t;x_0)$ be a solution of $\dot{x}(t) = f(t,x)$. If the perturbation
$y(t)$ is such that $x(t)+y(t)$ is also a solution of the same equation,
then $y(t)$ satisfies
\begin{equation*}
 \dot{y}(t) = f(t, y+x(t)) - f(t,x(t)) = F(t,y),
\end{equation*}
where clearly $F(t,0)\equiv 0$. The solution $x(t)$ is stable in the
sense of Lyapunov if for every $\epsilon > 0$ there exists a $\delta >
0$ such that $\norm{y(0)} < \delta$ implies $\norm{y(t)} < \epsilon$
for $t\geq 0$.  Now (see \cite{SC}, \cite{CL}), $x(t)$ is stable
if and only if the zero solution of $\dot{y}(t) = F(t,y)$ is
stable. Thus it is enough to look at zero solutions of systems of the
form $\dot{y}(t) = F(t,y)$, with $F(t,0) \equiv 0$.

Lyapunov pointed out an other possible simplification. Since $F(t,0)\equiv 0$,
if $F(t,y)$ is assumed to be continuously differentiable in $y$,
$F(t,y) = A(t)y + o(\norm{y})$ as $y\rightarrow 0$. Here,
\begin{align}
A(t) &= \frac{\partial F(t,y)}{\partial y} \Bigl\lvert_{y=0} \nonumber\\
     &= \frac{\partial f(t,x)}{\partial x} \Bigl\lvert_{x = x(t;x_0)}.
\label{eqn-fapprox}
\end{align}
Lyapunov argued that perhaps stability properties of the zero solution
of the linear first approximation $\dot{y}(t) = A(t)y$ might imply
stability properties of the zero solution of the nonlinear system
$\dot{y}(t) = F(t,y)$.  This program has been carried out by Lyapunov,
E. Cotton (1911), O. Perron (1928), I.G. Petrovskii (1934), R. Bellman (1953),
and others.
 
The other stream of research, which also originated with Lyapunov,
uses Lyapunov functions $V(t,y)$. If every solution $y(t; y_0)$ of
$\dot{y}(t) = F(t,y)$, $t\geq 0$, had the property that $\norm{y(t;y_0)}$
decreases as $t$ increases, stability of the zero solution would be
easy to infer. However, this property does not hold even for stable,
linear systems of the form $\dot{y}(t) = Ay$, where $A$ is a constant
matrix. Lyapunov's idea was to use a non-negative, real-valued
function $V(t,y)$, instead of the norm, such that $V(t,y(t))$
decreases as $t$ increases. This $V(t,y)$ is  required to be
related to $\norm{y}$ in a uniform way for $t\geq 0$; precise details
depend upon the concept of stability. After Lyapunov, this line was greatly
developed by Soviet researchers, including Persidskii, Malkin, Krasovskii,
and others, beginning in the 1920s. It was taken up by researchers in the
west, including J.L. Massera, T. Yoshizawa, and others, in the 1950s.

Parts of Sansone and Conti \cite{SC} and Hale \cite{Hale} are
excellent introductions to the theory of stability of ODEs. There are
numerous advanced works; of them we refer to Bellman \cite{Bellman},
Malkin \cite{Malkin}, and Yoshizawa \cite{Yoshizawa}.

Obtaining detailed stability information about solutions of ODEs can
be far from trivial. For this reason, applying the theory in Sections
5 and 6 to concrete examples is not a simple matter. In Section 7, we
give three applications to dynamical systems. We derive constant or
linear upper bounds on $E(t)$ for trajectories falling into stable,
hyperbolic fixed points, or into stable, hyperbolic cycles, or into a
normally hyperbolic and contracting manifold, with the flow on the
manifold being quasiperiodic.  These three applications involve the
Hartman-Grobman theorem, convergence in phase results for stable
cycles, and results of Pugh, Shub, and others about normally
hyperbolic flows, respectively. The first of these applications was
covered by Stuart and Humphries \cite{SH}. The second application is a
significant extension of a result due to Cano and Sanz-Serna
\cite{CanoS}.  The third application appears to be entirely new.

Let us mention the similarity of our analysis to the asymptotic analysis
of global error of Henrici \cite{Henrici1} \cite{Henrici2} and several
following researchers, including Gragg \cite{Gragg}.
Theorem \ref{thm-5-1}, in particular, is implicit in Henrici's work.
Beginning with Dahlquist \cite{Dahlquist}, it has been known that 
the use of one-sided Lipshitz constants can sometimes give meaningful
bounds on the global error. We show in Section 8 that the use of
one-sided Lipshitz constants fits naturally into our framework. 
Section 8 also considers variable time stepping, multistep methods,
and other issues.

To summarize briefly, this paper consists of three parts. The first 
part, Sections 2, 3, and 4, derives a conditioning function $E(t)$ and
associates it with the global error in numerically approximating the
solution of an ordinary differential equation. The second part, Sections
5 and 6, relates $E(t)$ to the stability theory of ordinary differential
equations. The third part, Section 7, applies this theory to dynamical
systems.

%% file: s1.ltx
\section{A Model for Discretization Errors}

The model of discretization error which we now present is close to the
way discretization errors are made by single step methods with
constant step sizes. Stuart and Humphries \cite{SH} model
discretization errors of single step methods in a similar manner.
 
\begin{figure}
\centerline{
\psfig{figure=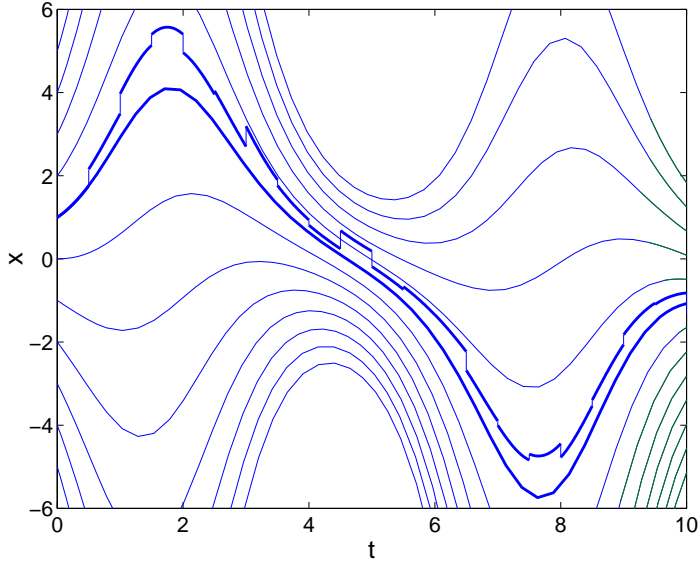,height=3in}}
\caption[xyz]{ The thick lines show an exact solution  of the equation
$\dot{x}(t) = \sin t + (\cos t) x$
and an approximation to it 
with $h=0.5$ and $\alpha(h) = 2h$ for.}
\label{fig-1}
\end{figure}

Let $\alpha(h)$ be a continuous, strictly increasing function of $h$ for 
$h\geq 0$. Assume also that $\alpha(0) = 0$. Then an approximation
$\tilde{x}_\alpha(t; x_0; h)$ to $x(t;x_0)$ is defined as follows:
\begin{align}
\tilde{x}_\alpha(0; x_0; h) &= x_0 \nonumber\\
\tilde{x}_\alpha(nh; x_0; h) &= x(nh; (n-1)h, 
\tilde{x}_\alpha((n-1)h; x_0; h))+h\alpha(h)v_n\quad n\geq 1,
\label{eqn-approx1}
\end{align}
where $v_n\in R^d$ can be any vector with $\norm{v_n} \leq 1$. In words,
the approximate solution at $t=nh$, $n\geq 1$, is obtained by exactly
propagating the point $\tilde{x}_\alpha((n-1)h;x_0;h)$ at $t=(n-1)h$
under $\dot{x}(t) = f(t,x)$ till $t=nh$, and then adding the discretization
error or discontinuity $h\alpha(h)v_n$, where $\norm{v_n}\leq 1$. For
$(n-1)h \leq t < nh$,
\begin{equation}
\tilde{x}_\alpha(t;x_0; h) = x(t; (n-1)h, \tilde{x}_\alpha((n-1)h; x_0; h)).
\label{eqn-approx2}
\end{equation}
Since $v_n$ can be any vector with $\norm{v_n} \leq 1$, this actually
defines a whole family of approximate solutions which we  denote
by $\tilde{X}_\alpha(x_0; h)$. The exact solution $x(t)$ is the only
member of this family which is continuous. Figure \ref{fig-1} gives an example
of an approximate solution.

We now address how single step numerical methods are related to 
approximations from the family $\tilde{X}(x_0; h)$. In our description
of single step methods, the discretization error at every step was
taken to be $K_ih^{r+1}v_i$. We now assume that the $K_i$ are bounded
by a constant $K$ which does not depend upon $h$ or $i$. This can be
proven in some circumstances; see \cite{SH}. Besides, if there is no
such $K$, the numerical method will not in practice behave as if it
were of order $r$.  Thus when the order of accuracy of the numerical
method for approximating $x(t;x_0)$ is $r$, we can take $\alpha(h) =
Kh^{r}$, and there will always be an approximation in the family
$\tilde{X}_\alpha(x_0;h)$ which is the same as the trajectory of the
numerical method.

The notations we have established so far will be adhered to in the rest
of the paper. For the system, $\dot{x}(t) = f(t,x)$, $f(t,0)$ is not
necessarily zero. The solution with $x(t_0) = x_0$ is denoted by 
$x(t; t_0, x_0)$, and by $x(t;x_0)$ when $t_0=0$. The approximations to
$x(t; x_0)$ obtained as in \eqref{eqn-approx1} and \eqref{eqn-approx2}
are denoted by $\tilde{x}_\alpha(t; x_0; h)$. The family of approximations
is $\tilde{X}(x_0; h)$. The same notations apply for $\dot{y}(t) = F(t,y)$
with the $x$s changed to $y$s. But here $F(t,0) \equiv 0$, and 
$\tilde{y}_\alpha(t;h)$, $t\geq 0$, is an approximation to the zero solution,
and $\tilde{Y}_\alpha(h)$ is a collection of those approximations. When
we speak of a linear first approximation, it is always obtained as in
\eqref{eqn-fapprox}. Sometimes we omit the word linear. This same equation
is sometimes called the equation of first variation or simply
the linearization. Let us note that
when we speak of stability, it is always stability in the sense of Lyapunov
and his followers; in particular, it is not numerical stability in the sense
of Dahlquist.

%% file: s2.ltx
\section{Definition of $E_\alpha(t)$}

Let us recall that $\alpha(h)$ is assumed to be a strictly increasing,
continuous function of $h$ for $h\geq 0$ which is zero at zero; for
example, $\alpha(h)$ can be $Kh^r$ for a positive integer $r$.  Assume
the functions $f(t,x)$ and $F(t,x)$ to be defined for $0\leq t <
\infty$ and $x\in R^d$, to be continuous in $t$, and locally Lipshitz
in $x$. Besides, $F(t,0)\equiv 0$. Assume the solution $x(t;x_0)$ of
the initial value problem $\dot{x}(t) = f(t,x)$, $x(0) = x_0$, to
be continuable till $t=\infty$. Only some of these assumptions are
restated in theorems that follow.

The global error $e_\alpha(t;x_0;h)$ is defined as follows:
\begin{equation*}
e_\alpha(t;x_0;h) = \sup_{\tilde{x}_\alpha\in \tilde{X}_\alpha}
\norm{\tilde{x}_\alpha(t;x_0;h) - x(t; x_0)}.
\end{equation*}
We later use this to define $E_\alpha(t)$.



Instead of saying that $\tilde{x}_\alpha(t;x_0;h)$ exists for $0\leq
t\leq T$, we say that $\tilde{x}_\alpha(t;x_0;h)$ can be continued
till $t=T$. We say that every approximation
$\tilde{x}_\alpha(t;x_0;h)$ can be continued till $T$ if every
$\tilde{x}_\alpha(t;x_0;h)$ exists for $0\leq t\leq T$ with any
allowed choice of discontinuities at $t=kh$. Lemma \ref{lem-2-1}
introduces $h_0(T,r)$ and $L(T,r)$.

\begin{lem}
Assume, as usual, that $f(t,x)$ is continuous in $t$, $t\geq 0$, and
locally Lipshitz in $x$, $x\in R^d$. Let $x(t;x_0)$, $t\geq 0$, be the
unique solution of the initial value problem $\dot{x}(t) = f(t,x)$,
$x(0) = x_0$. Then,
\begin{description}
\item[(i)]  There exists a constant $L(T,r)>0$ such that 
\begin{equation*}
\norm{f(t,x_1) - f(t,x_2)} \leq L(t,r) \norm{x_1 - x_2}
\end{equation*}
for $0\leq t\leq T$ and $\norm{x_i - x(t;x_0)} \leq r$, for any
$T>0$, $r>0$, and $i=1,2$.

\item[(ii)] There exists a constant $h_0(T,r)$ such that for $0< h < h_0(T,r)$
every approximation $\tilde{x}_\alpha(t;x_0;h)$ can be continued till 
$t=T$, and satisfies $\norm{\tilde{x}_\alpha(t;x_0; h)-x(t;x_0} < r$,
for $0\leq t\leq T$.

\end{description}
Further, for $0<h<h_0(T,r)$, $e_\alpha(t; x_0; h) \leq t e^{L(T,r)t}\alpha(h)$.
\label{lem-2-1}
\end{lem}
\begin{proof}
Similar to proof of Theorem 3.4.6 in Stuart and Humphries \cite{SH}. Similar
estimates using the Lipshitz constant are found in \cite{HNW} and other
places. 
\end{proof}

For the zero solution of $\dot{y}(t) = F(t,y)$, $y(0)=0$, where
$F(t,0)\equiv 0$, $e_\alpha(t;h)$ is defined as follows:
\begin{equation}
 e_\alpha(t;h) = \sup_{\tilde{y}_\alpha \in \tilde{Y}_\alpha}
\norm{\tilde{y}_\alpha(t;h)},
\label{eqn-defn-e}
\end{equation}
where $\tilde{Y}_\alpha(h)$ is the family of approximations to the zero
solution  with time step $h$.

\begin{prop}
Assume $h < h_0(T,r)$ for $r>0$, $T>0$, and that $0\leq t \leq T$.  The
$e_\alpha(t;x_0;h)$ of the solution $x(t;x_0)$ of $\dot{x}(t) =
f(t,x)$, $x(0) = x_0$, and the $e_\alpha(t;h)$ of the zero solution of
$\dot{y}(t) = F(t,y)$, $y(0)=0$, where $F(t,y) = f(t,y+x(t;x_0)) -
f(t,x(t;x_0))$, are the same.
\label{prop-2-1}
\end{prop}

\begin{proof}
It is enough to show that members $\tilde{x}_\alpha(t;x_0;h)$ of
$\tilde{X}_\alpha(x_0;h)$ and members $\tilde{y}_\alpha(t;h)$ of
$\tilde{Y}_\alpha(h)$ can be matched so that
$\tilde{x}_\alpha(t;x_0;h) = x(t;x_0) + \tilde{y}_\alpha(t;h)$.

It is well known (see \cite{CL}) that
\begin{equation*}
x(t+\tau; t; x(t;x_0)+\delta) = x(t+\tau;x_0) + y(t+\tau;t,\delta)
\end{equation*}
for $\tau\geq 0$. Hence, if $\tilde{x}_\alpha(kh; x_0; h) =
x(kh;x_0) + \tilde{y}_\alpha(kh;h)$ then $\tilde{x}_\alpha(t; x_0; h) =
x(t;x_0) + \tilde{y}_\alpha(t;h)$ for $kh < t < (k+1)h$. Further, the
discontinuities of $\tilde{x}_\alpha(t;x_0;h)$ and $\tilde{y}_\alpha(t;h)$
at $t=kh$, where $k$ is a positive integer, can be exactly matched. Therefore,
the approximate solutions can be matched as desired.
\end{proof}

The $E_\alpha(t)$ which corresponds to the zero solution of
$\dot{y}(t) = F(t,y)$, $y(0) = 0$, $t\geq 0$, is defined as
follows:
\begin{equation}
 E_\alpha(t) = \limsup_{h\rightarrow 0} \frac{e_\alpha(t;h)}{\alpha(h)}.
\label{eqn-defn-E}
\end{equation}
By Lemma \ref{lem-2-1}, $E_\alpha(t) \leq t e^{Lt}$, where $L=L(T,r)$
for some $T>t$ and $r>0$. For a nonzero solution $x(t;x_0)$ of
$\dot{x}(t) = f(t,x)$, $x(0) = x_0$, $t\geq 0$, $E_\alpha(t)$ is
defined to be the same as the $E_\alpha(t)$ for the zero solution of
$\dot{y}(t) = f(t,y+x(t;x_0)) - f(t, x(t;x_0))$, $y(0)=0$, $t\geq 0$.

As in the stability theory of ODEs, we can and do confine ourselves to
an analysis of the zero solution of $\dot{y}(t) = F(t,y)$ without any
loss of generality.  From here on, assume $L(T,r)$ to be such that
$$\norm{F(t,y_1) - F(t,y_2)} \leq L(t,r) \norm{y_1-y_2}$$ for $0\leq
t\leq T$, $\norm{y_1} \leq r$, $\norm{y_2}\leq r$, where $T>0$ and
$r>0$. Also, assume $h_0(T,r)$ to be such that for $0<h < h_0(T,r)$
every approximation $\tilde{y}_\alpha(t;h)$ can be continued till
$t=T$ and satisfies $\norm{\tilde{y}_\alpha(t;h)} < r$.

%% file: s3.ltx
\section{Properties of $E_\alpha(t)$}

This section is devoted to properties of $E_\alpha(t)$.
We show that $E_\alpha(t)$ is continuous (Theorem 
\ref{thm-3-1}) and independent of the choice of $\alpha(h)$ (Theorem
\ref{thm-3-2}) so that it can be written as $E(t)$. Theorem \ref{thm-3-3}
and its corollary relate $E(t)$ to the accumulation of global error
when $x(t;x_0)$ is approximated.

In the proofs in this section, we repeatedly use Theorem 10.1 of \cite{HNW}.
That theorem allows us to bound the divergence of two solutions of 
a differential equation using a Lipshitz constant. 

\newcommand{\ea}[1]{e_\alpha(#1;h)}
\newcommand{\ya}[1]{\tilde{y}_\alpha(#1;h)}

A sequence of inequalities are often combined in a way we now
describe.  Let $e_0=0$, and $e_i \leq f e_{i-1} + r_{i-1}$ for $1\leq
i\leq n$.  Then, we have,
\begin{equation*}
e_n \leq (f^{n-1}r_0 + f^{n-2}r_1 +\cdots+ r_{n-1}),
\end{equation*}
assuming $f\geq 0$. In fact, most of the time $r_i$ will all be
equal. In that situation, $e_n \leq r_0(1+f+\cdots+f^{n-1})$. And when
$f\geq 1$, we can get $e_n \leq f^{n-1}(r_0+\cdots+r_{n-1})$, or $e_n
\leq nf^{n-1}r_0$ if the $r_i$ are all equal.

\begin{lem}
Let $h < h_0(t+s,r)$, $r>0$. Then
\begin{equation*}
\ea{t+s} \leq e^{Ls} \ea{t} + (s+h)e^{Ls}\alpha(h),
\end{equation*}
where $L=L(t+s,r)$ and $s\geq 0$.
\label{lem-3-1}
\end{lem} 
\begin{proof}
Let $\ya{t}$ be an approximation to the zero solution of $\dot{y}(t) =
F(t,y)$, $y(0)=0$. 

Let $t+h_l$ be the first multiple of $h$ after $t$, $t+s-h_r$ the last
multiple before $t+s$, and let $t+s-h_r = t+h_l + nh$. Using Theorem 10.1
of \cite{HNW} and taking into account the discontinuities at multiples
of $h$, we have
\begin{align*}
\norm{\ya{t+h_l}} &\leq e^{Lh_l}\norm{\ya{t}} + h\alpha(h), \\
\norm{\ya{t+h_l+kh}} &\leq e^{Lh}\norm{\ya{t+h_l+(k-1)h}}+h\alpha(h),
\quad \text{$1\leq k\leq n$}, \\
\norm{\ya{t+s}} &\leq e^{Lh_r}\norm{\ya{t+nh}}.
\end{align*}
Combining these inequalities, we get
\begin{equation*}
\norm{\ya{t+s}} \leq e^{Ls} \norm{\ya{t}} + (s+h) e^{Ls} \alpha(h).
\end{equation*}
The proof is easy to complete using \eqref{eqn-defn-e}. 
\end{proof}

\begin{lem}
With the same assumptions about $h$ and $L$ as in the previous lemma,
\begin{equation*}
\ea{t+s} \geq e^{-Ls} \ea{t},
\end{equation*}
where $s\geq 0$.
\label{lem-3-2}
\end{lem}
\begin{proof}
Let $\ya{\tau}$, $0\leq \tau\leq t$, be an approximate solution. Continue 
it from $\tau=t$ to $\tau=t+s$ exactly as a solution of $\dot{y}(t)
= F(t,y)$. Then,
\begin{equation*}
\norm{\ya{t+s}} \geq e^{-Ls} \norm{\ya{t}}.
\end{equation*}
The proof is now easy to complete using \eqref{eqn-defn-e}.
\end{proof}

\begin{thm}
The $E_\alpha(t)$ of the zero solution of
$\dot{y}(t) = F(t,y)$, $y(0)=0$, is continuous for $t\geq 0$. 
\label{thm-3-1}
\end{thm}
\begin{proof}
From Lemmas \ref{lem-3-1} and \ref{lem-3-2}, we get
\begin{equation*}
\ea{t-\delta} e^{-L\delta} \leq \ea{t} \leq e^{L\delta}\ea{t-\delta} 
+ e^{L\delta}(\delta +h)\alpha(h),
\end{equation*}
when $0 < t-\delta < t$, and when $t < t+\delta$,
\begin{equation*}
(\ea{t+\delta} - e^{L\delta}(\delta+h)\alpha(h))e^{-L\delta} \leq
\ea{t} \leq e^{L\delta} \ea{t+\delta}.
\end{equation*}
Divide these two inequalities by $\alpha(h)$ and use \eqref{eqn-defn-E}
to deduce continuity of $E_\alpha(t)$.
\end{proof}

\newcommand{\eA}[1]{e_\alpha(#1;h_\alpha)}
\newcommand{\yA}[1]{\tilde{y}_\alpha(#1;h_\alpha)}
\newcommand{\eB}[1]{e_\beta(#1;h_\beta)}
\newcommand{\yB}[1]{\tilde{y}_\beta(#1;h_\beta)}

\begin{thm}
The $E_\alpha(t)$ of the zero solution of
$\dot{y}(t) = F(t,y)$, $y(0)=0$, is the same for any $\alpha(h)$, with
$\alpha(0) = 0$ and $\alpha(h)$ continuous and strictly increasing for
$h\geq 0$.
\label{thm-3-2}
\end{thm}
\begin{proof}
Let $\beta(h)$ be another function like $\alpha(h)$, which is
continuously increasing for $h\geq 0$ and zero at zero. We will show
that $E_\beta(t) = E_\alpha(t)$.

Take $L=L(t,r)$ for some $r>0$. Take $h_\alpha$ and $h_\beta$ small enough
that the approximate solutions in the families $\tilde{Y}_\alpha(h)$
and $\tilde{Y}_\beta(h)$ are not only continuable till $t$ but stay within
a radius $r$ of zero. For $\epsilon$ small enough, there exist unique
$h_\alpha$ and $h_\beta$ such that $\alpha(h_\alpha) = \beta(h_\beta)
= \epsilon$. Assume this relationship between $h_\alpha$, $h_\beta$,
and $\epsilon$. Then $h_\alpha\rightarrow 0$ and $h_\beta\rightarrow 0$
as $\epsilon\rightarrow 0$. 

It is enough to show that 
\begin{equation*}
\lim_{\epsilon\rightarrow 0} \frac{\abs{\eA{t} - \eB{t}}}{\epsilon} = 0.
\end{equation*}
In fact, it is enough to show that for every approximate solution 
$\yA{t}$ there is an approximate solution $\yB{t}$ such that
\begin{equation}
\norm{\yA{t} - \yB{t}} \leq \epsilon f(t) g(\epsilon),
\label{eqn-s3-1}
\end{equation}
where $f(t)$ is finite and $g(\epsilon)\rightarrow 0$ as $\epsilon\rightarrow
0$. It will also have to be shown that a given $\yB{t}$ can be approximated
by some $\yA{t}$ in the same manner; but the proof of this is gotten by
transposing $\alpha$ and $\beta$. 

Given a $\yA{t}$, we construct a $\yB{t}$ as follows. $\yB{\tau}$ is
determined by the differential equation $\dot{y}(t) = F(t,y)$, except
when $\tau = ph_\beta$, $p$ being a positive integer. At points 
$\tau = ph_\beta$, introduce a discontinuity of magnitude smaller than
or equal to $h_\beta\beta(h_\beta)$ in such a way as to get as close to
$\yA{ph_\beta}$ as possible.

Let $h_\beta= (m+f)h_\alpha$, $m\geq 0$ being an integer, and $0\leq f<1$.
Assume that the number of integer multiples of $h_\alpha$ in 
$(ph_\beta,(p+1)h_\beta]$ is $k_p$. Let $e_p = \norm{\yA{ph_\beta}
-\yB{ph_\beta}}$.

Then $e_0 = 0$ and
\begin{align*}
e_{p+1} &\leq e^{h_\beta L} e_p + k_p h_\alpha \alpha(h_\alpha)e^{h_\beta L}
	- h_\beta \beta(h_\beta) \\ 
	&= e^{h_\beta L}e_p + k_p h_\alpha\epsilon e^{h_\beta L} -
  	h_\beta \epsilon.
\end{align*}
The first two terms can be derived from Theorem 10.1 of \cite{HNW}. The
third term is because of the discontinuity introduced into $\yB{\tau}$ 
at $\tau=ph_\beta$. If this bound on $e_{p+1}$ is negative, we can take
$e_{p+1}=0$. If $n =\floor{\frac{t}{h_\beta}}$, and there are $k_n$ multiples
of $h_\alpha$ in $(nh_\beta,t]$, then
\begin{equation*}
\norm{\yA{t} - \yB{t}} \leq e^{Lh_r}e_n + k_n h_\alpha\alpha(h_\alpha)e^{Lh_r},
\end{equation*}
where $h_r = t-nh_\beta$.

Let $p$ be the highest among $\{0,1,\ldots,n\}$ such that $e_p=0$. Combining
the inequalities for $e_i$, $p\leq i\leq n$, we get
\begin{equation*}
e_n \leq e^{L(n-p)h_\beta}(h_\alpha\epsilon e^{h_\beta L}(k_p+\cdots+k_{n-1})
 	-(n-p)h_\beta\epsilon).
\end{equation*}
But $k_p+\cdots+k_{n-1}$ being the number of multiples of $h_\alpha$ in
$(ph_\beta, nh_\beta]$ is bounded above by $(n-p)h_\beta/h_\alpha+1$.
Let $t'=t-h_r$. Then,
\begin{align*}
e_n &\leq e^{Lt'}((n-p)h_\beta e^{h_\beta L}\epsilon + h_\alpha\epsilon
 e^{h_\beta L} - (n-p)h_\beta\epsilon) \\
    &\leq \max(t'e^{Lt'}, e^{Lt'}) ((e^{h_\beta L}-1) + 
    h_\alpha e^{h_\beta L})\epsilon.
\end{align*}
Finally,
\begin{equation*}
\norm{\yA{t} - \yB{t}} \leq \max(e^{Lt},te^{Lt}) ((e^{h_\beta L}-1)
+ h_\alpha e^{h_\beta L} + k_n h_\alpha)\epsilon.
\end{equation*}
Since $k_n$ is the number of multiples of $h_\alpha$ in $(t-h_r,t]$,
$k_n h_\alpha <h_\beta+h_\alpha$. Thus the bound above satisfies all
conditions required of \eqref{eqn-s3-1}, and the proof is complete.
\end{proof} 

\newcommand{\yt}[1]{\tilde{y}(#1;h)}
From here on, we drop the subscript $\alpha$ from $E_\alpha(t)$,
$\tilde{y}_\alpha(t;h)$, $\tilde{Y}_\alpha(h)$,
$\tilde{x}_\alpha(t;x_0,h)$, $\tilde{X}_\alpha(x_0;h)$, and $e_\alpha(t;h)$.

\begin{thm}
Let $E(t)$ be associated with the zero solution of $\dot{y}(t)=F(t,y)$,
$y(0)=0$. Given $T>0$ and $\epsilon>0$, there exists $h_0>0$ such that
$h < h_0$ implies 
\begin{equation*}
\sup_{\tilde{y}\in\tilde{Y}}\norm{\yt{t}} = e(t;h) 
\leq (E(t)+\epsilon)\alpha(h)
\end{equation*}
for $0\leq t\leq T$.

Further, if $\theta< E(t)$ there exists an $\epsilon > 0$ such that
$e(t;h)>(\theta+\epsilon)\alpha(h)$ for arbitrarily small $h$.
\label{thm-3-3}
\end{thm}

\begin{proof}
Let
$\epsilon(t;h) = e(t;h)/\alpha(h) - E(t)$.

Take $h<h_0$, where $h_0$ to begin with is smaller than $h_0(T,r)$ for
some $r>0$, and $L=L(T,r)$.  By definition of $E(t)$
\eqref{eqn-defn-E},
\begin{equation}
\limsup_{h\rightarrow 0} \epsilon(t,h) = 0
\label{eqn-s3-2}
\end{equation}
for $0\leq t\leq T$. 

To prove the first part, 
given $\epsilon > 0$ we find an $h_0>0$ such that $h<h_0$ implies 
$\epsilon(t;h)< \epsilon$ for $0\leq t\leq T$.
In fact, it is enough to show that there exists an $h_t>0$,
which depends on $t$, and an open neighbourhood of $t$ in $[0,T]$
such that for $\tau$ in that neighbourhood and $h<h_t$, $\epsilon(\tau;h)
< \epsilon$. If the neighbourhoods of $t_1,\ldots,t_n$ give a finite
cover of the compact interval $[0,T]$, we can take $h_0$ to be the
minimum of $h_{t_1},\ldots,h_{t_n}$.

Clearly,
\begin{align*}
\epsilon(t+\delta t;h) &= \epsilon(t;h) +
\frac{e(t+\delta t;h) - e(t;h)}{\alpha(h)} -
(E(t+\delta t) - E(t))\\
&\leq \epsilon(t;h) + \Abs{\frac{e(t+\delta t;h) - e(t;h)}{\alpha(h)}}
+ \abs{E(t+\delta t) - E(t)}.
\end{align*}
By \eqref{eqn-s3-2}, $\epsilon(t;h) < \epsilon_1$ for $h$ small enough.
By continuity of $E(t)$, $\abs{E(t+\delta t) - E(t)}<\epsilon_2$ for
$\abs{\delta t}$ small enough.

Using Lemmas \ref{lem-3-1} and \ref{lem-3-2}, and the bound
$tE^{Lt}\alpha(h)$ on $e(t;h)$ from Lemma \ref{lem-2-1}, we get
\begin{equation*}
\frac{\abs{e(t+\delta t;h) - e(t;h)}}{\alpha(h)}\leq
\max\Bigl(te^{Lt}(1-e^{-L\delta t}), te^{Lt}(e^{L\delta t} - 1)
+ (\delta t + h) e^{L\delta t}\Bigr).
\end{equation*}
for $\delta t >0$, and for $\delta t < 0$,
\begin{equation*}
\frac{\abs{e(t+\delta t;h) - e(t;h)}}{\alpha(h)} \leq
\max\Bigl(te^{Lt}(e^{-L\delta t} -1), te^{Lt}(1-e^{L\delta t})+
     (\delta t +h)\Bigr).
\end{equation*}
Therefore, we can choose $\abs{\delta t}$ and $h$ small enough to ensure
$\abs{e(t+\delta;h) - e(t;h)}/\alpha(h) < \epsilon_3$. We need only
take $\epsilon_1 + \epsilon_2 + \epsilon_3 < \epsilon$ to find the desired
$h_t$ and the neighbourhood of $t$. 

The second half of the theorem is immediate from the definition of $E(t)$.
\end{proof}

\begin{cor}
Let $E(t)$ be associated with the solution $x(t;x_0)$ of $\dot{x}(t)
= f(t,x)$, $x(0)=x_0$, where $f(t,x)$ is continuous in $t$ and locally
Lipshitz in $x$ in an open neighbourhood of the solution $(t,x(t;x_0))$,
$t\geq 0$. Given $T>0$ and $\epsilon > 0$, there exists an $h_0>0$ 
such that 
\begin{equation*}
\norm{\tilde{x}(t;x_0;h)-x(t;x_0)} < (E(t)+\epsilon)\alpha(h)
\end{equation*}
for $0\leq t\leq T$, $h<h_0$, and any approximation $\tilde{x}(t;x_0;h)$
in the family $\tilde{X}(x_0;h)$.
\label{cor-3-1}
\end{cor}

\begin{figure}
\centerline{
\psfig{figure=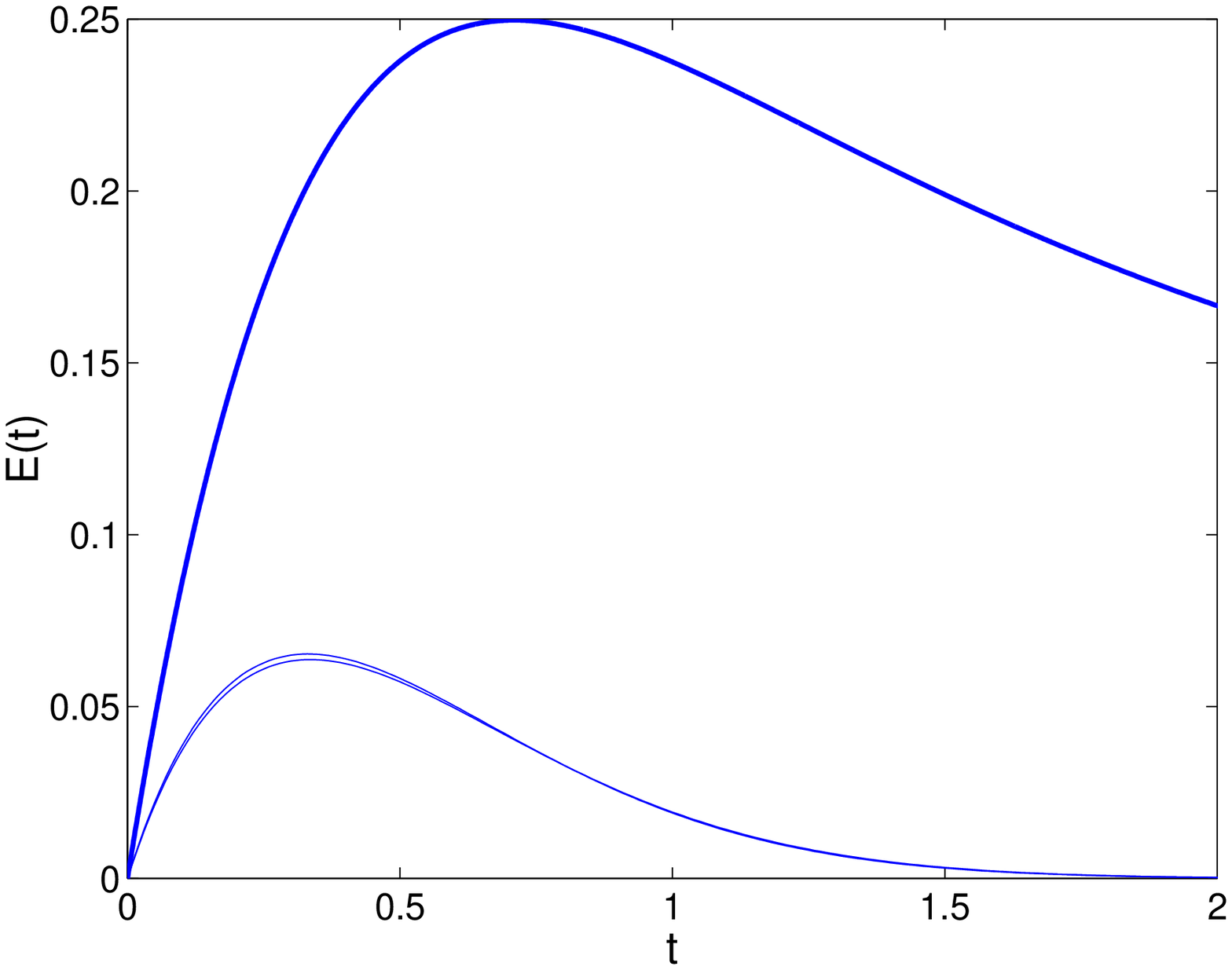,height=2.5in}
\psfig{figure=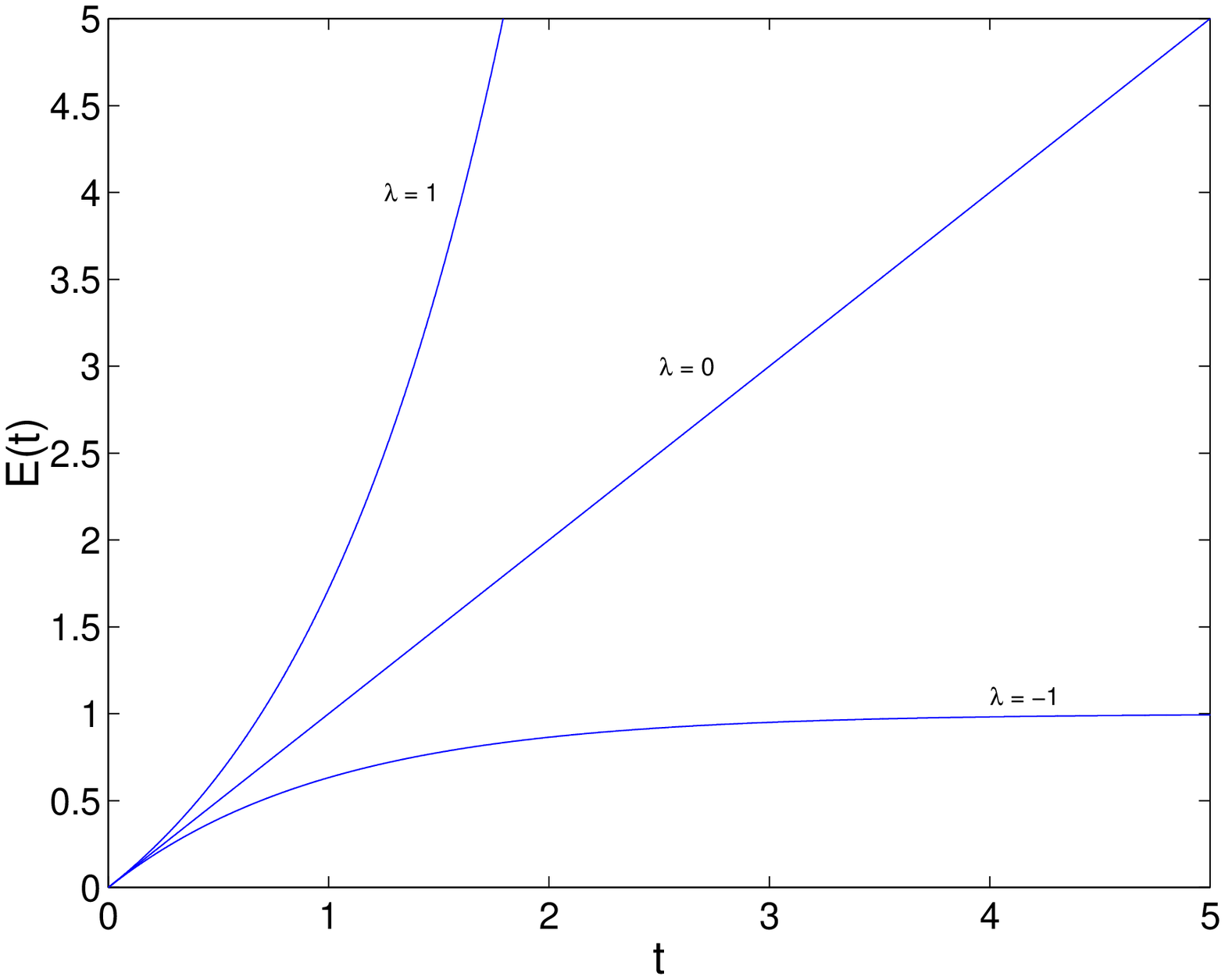,height=2.5in}}
\centerline{
\hspace{.2in} {\bf (a)} \hspace{2.8in} {\bf (b)}}
\caption[xyz] { (a) The thick line  is the $E(t)$ for $\dot{y}(t)
= -(2(t+1) + (t+1)^{-1})y$, $y(0)=1$. The two thin lines below it are
global errors of forward and backward Euler divided by $8 h$; here
$h=.01$. (b) $E(t)$ for $\dot{y}(t) = \lambda y$. The exact $E(t)$ are
obtained using Theorem \ref{thm-4-1}.}
\label{fig-2}
\end{figure}

So far, we have defined the notion of an approximate solution (see
Figure \ref{fig-1}), and an $E(t)$ for every solution of an
ODE. Theorem \ref{thm-3-3} and its Corollary \ref{cor-3-1} relate
$E(t)$ to the accumulation of global error. In the upper bound
$(E(t)+\epsilon)\alpha(h)$ for the global error, the details of the
numerical method have been pushed into $\alpha(h)$; the conditioning
function $E(t)$ is something intrinsic to the solution
$x(t;x_0)$. Figure \ref{fig-2} shows $E(t)$ for some simple examples.

As shown in Figure \ref{fig-2}a, the bound on global error given by
$E(t)$ may be pessimistic in practice. Our model of approximations
allows arbitrary discontinuities whose norms are bounded
by $h\alpha(h)$, while for any given numerical method
the discretization errors are fixed. Our analysis does not say how
small an $h$ is good enough for the bound on the global error to be
governed by $E(t)$ as in Theorem \ref{thm-3-3}.  Will an $h$ in a
practical computation be small enough? We answer this in Henrici's
words \cite{Henrici2}: {\it usually if a stepsize is small enough to
yield an accurate solution, it is also small enough that an asymptotic
formula gives a correct indication of the size of the error.}
Numerical experiments in \cite{CalvoS} and \cite{CanoS} lend ample
support to this statement.

%% file: s4.ltx
\section{$E(t)$ for Linear Systems}

Our investigation of the relationship between $E(t)$ and stability 
properties of the exact solution begins with the linear system 
$\dot{y}(t) = A(t)y$, $y(0)=0$. The relationship is not as simple as
one might wish. There are both asymptotically stable examples with
exponentially increasing $E(t)$ and unstable examples with linearly
bounded $E(t)$. However, Theorems \ref{thm-4-3} and \ref{thm-4-4}
give conditions for $E(t)$ to be bounded by a constant or to be 
linearly bounded. Linear systems are an important class of problems by
themselves. What is more, they can be used for understanding the $E(t)$
of nonlinear systems. We will see this in Theorem \ref{thm-5-1}
and in Section 6.2.

In $\dot{y}(t) = A(t)y$, 
the $d\times d$ matrix $A(t)$ is assumed to be continuous for $t\geq
0$.  For such linear systems, it is easy to show that $y(t)$ is a
linear function of $y(0)$ for $t\geq 0$ \cite{SC}. Thus
$y(t)=Y(t)y(0)$ for $Y(t)\in R^{d,d}$, where $Y(t)$ is called the
principal fundamental matrix of the linear system. The matrix $Y(t)$
itself is continuous in $t$ and always nonsingular. Moreover, $y(t) =
Y(t)\inv{Y}(s)y(s)$ for $t\geq s$. Although it does not appear to be
directly related, let us mention the beautiful derivation of the
Magnus series for $Y(t)$ by Iserles and N{\o}rsett \cite{IN}.

For scalar linear systems $\dot{y}(t) = a(t)y$, $a(t)\in R$, we have the 
following theorem. 

\begin{thm}
The $E(t)$ for the zero solution of $\dot{y}(t)=a(t)y$, $y(0)=0$, is
given by 
\begin{equation*}
E(t) = e^{g(t)} \int_{0}^{t} e^{-g(s)} ds,
\end{equation*}
where
\begin{equation*}
g(t) = \int_0^{t}a(\tau) d\tau.
\end{equation*}
\label{thm-4-1}
\end{thm}
\begin{proof}
The fundamental matrix, which is  scalar in this situation, is given by
$Y(t)=e^{g(t)}$. Since $Y(t)$ is always positive, the optimal choice
of $v(s)$ in Theorem \ref{thm-4-2} is $v(s)\equiv 1$.
\end{proof}

 Let us now consider the concepts of stability put forward by Lyapunov
\cite{Lyapunov}. The definitions will be stated for nonzero solutions
$x(t;x_0)$.

\begin{defn}
The solution $x(t;x_0)$ of $\dot{x}(t)=f(t,x)$, $x(0)=x_0$, is {\it stable}
if given any $\epsilon>0$ there exists a $\delta > 0$ such that
$\norm{x_0' - x_0}<\delta$ implies $\norm{x(t;x_0')-x(t;x_0)} < \epsilon$  
for $t\geq 0$. In fact, stability implies that given $\epsilon > 0$,
there exists a $\delta(\tau)>0$ for every $\tau\geq 0$ such that
$\norm{x(\tau;x_0') - x(\tau; x_0)} < \delta(\tau)$ implies that
$\norm{x(t;x_0') - x(t;x_0)} < \epsilon$ for $t\geq \tau$. Let us
emphasize that $\delta(\tau)$ can depend on $\tau$.
\end{defn}

\begin{defn}
The solution $x(t;x_0)$ is {\it asymptotically stable} if given $\epsilon > 0$
there exists a $\delta(\tau)>0$ for every $\tau\geq 0$ such that
$\norm{x_\tau'-x(\tau;x_0)}<\delta(\tau)$ implies not only that 
$\norm{x(t;\tau,x_\tau') - x(t;x_0)} < \epsilon$ for $t\geq \tau$ 
but also that $\norm{x(t;\tau,x_\tau') - x(t;x_0)}\rightarrow 0$ as
$t\rightarrow\infty$. 
\end{defn}

Implicit in the definitions is an assumption about the existence of
solutions which begin near the solution $x(t;x_0)$. Obviously, 
asymptotic stability implies stability. For the scalar, linear problem
$\dot{y}(t) = a(t)y$, $y(0)=y_0$, a necessary and sufficient condition
for asymptotic stability is $g(t)\rightarrow -\infty$ as $t\rightarrow\infty$,
where $g(t) = \int_0^t a(s) ds$. However, the following examples show
that both these concepts of stability are insufficient for bounding
$E(t)$.

\begin{exmp}
We will show using the scalar, linear problem $\dot{y}(t) = a(t)y$ that
even when a solution is asymptotically stable the $E(t)$ associated with
it can increase at an arbitrarily high rate. Given a rate $r(t)$,
consider a continuously differentiable function $g(t)$, $t\geq 0$,
$g(0) = 0$, such that
\begin{enumerate}
\item $g(t) \leq -t$ for all $t\geq 0$,
\item $e^{g(k)}\int_{0}^{k}e^{-g(s)}ds > r(k)$ for $k=1,2,3,\ldots$
\end{enumerate}
For the linear system, take $a(t) = g'(t)$. The first condition
ensures asymptotic stability of the zero solution, and the second
condition implies $E(k) > r(k)$ for positive integers $k$.  Such a
$g(t)$ is easy to construct. Take $g(k)=-k$ for $k=0,1,2,\ldots$. For
$k-1 < t < k$, $k\geq 1$, define $g(t)$ so that $g(t) \leq -t$ and
\begin{equation*}
\int_{k-1}^{k} e^{-g(s)} ds \geq r(k) e^k.
\end{equation*}
This can be carried out for any continuous $r(t)$, for example $r(t)=e^t$.
\label{exmp-4-1}
\end{exmp}

\begin{exmp}
On the other hand, there are unstable solutions with linearly bounded
$E(t)$. Consider the scalar, linear system $\dot{y}(t) = \frac{\alpha}{t+1}
y$, $t\geq 0$. For this ODE, $y(t) = (1+t)^\alpha y(0)$ implying instability
of the zero solution for $\alpha > 0$. Yet, for $0<\alpha <1$
$E(t)$, which is $(1-\alpha)^{-1} (1+t)^\alpha((1+t)^{1-\alpha}-1)$, is
linearly bounded. For $\alpha=1$, $E(t)$ is $(1+t)\log(1+t)$.
\label{exmp-4-2}
\end{exmp}

\begin{exmp}
In Example \ref{exmp-4-1}, $\abs{a(t)}$ will be unbounded. Does asymptotic
stability of $\dot{y}(t) = a(t)y$, $t\geq 0$, imply a linear bound for
$E(t)$ if $\abs{a(t)}$ is bounded? The answer is no; $E(t)$ can still
increase exponentially in $t$. We sketch the construction of a $g(t)$
to show this. First take $g_1(t) = -t$ and $g_2(t) = -2t$. Take
$g(t) = g_2(t)$ for $0 \leq t \leq t_1$, and let $g(t)$ increase monotonically
till $g(t_2) = g_1(t_2)$ for $t_2 \geq t_1$, and then let $g(t)$ decrease
monotonically till $g(t_3) = g_2(t_3)$ for $t_3\geq t_2$. Repeat the same
construction from $t_3$ onwards with $t_4$, $t_5$, and $t_6$ in place of
$t_1$, $t_2$, $t_3$, and so on. The construction may be arranged so that
\begin{enumerate}
\item if $g(\tau) = g_1(\tau)$ then $g(t) = g_2(t)$ for 
$f_1 \tau \leq t \leq f_2 \tau$ for any fixed $0 < f_1 < f_2 < 1$,
\item $\abs{a(t)} = \abs{g'(t)}$ is bounded. 
\end{enumerate}
It is easy to check that $E(\tau) \geq e^{-\tau/2}(e^{2f_2\tau} -
e^{2f_1\tau})$, for $\tau$ such that $g(\tau) = g_1(\tau)$. Therefore,
for $f_2 > 1/2$, $E(t)$ increases exponentially. Let us note that the
linear system in this example has a negative Lyapunov exponent of
$-1$.
\label{exmp-4-3}
\end{exmp} 

Notions of stability needed for bounding $E(t)$ either by a constant or 
linearly are introduced after Theorem \ref{thm-4-2}. Theorem \ref{thm-4-2}
is comparable to Theorem 3.1 of Cano and Sanz-Serna \cite{CanoS}.

\begin{thm}
The $E(t)$ of the zero solution of $\dot{y}(t) = A(t)y$, $y(0)=0$,
is given by 
\begin{equation*}
 E(t) = \sup_{v(s)} \Norm{\int_0^t Y(t)\inv{Y}(s) v(s) ds},
\end{equation*}
where the supremum is over all continuous functions
$v(s):[0,t]\rightarrow R^d$ with $\norm{v(s)}\leq 1$ for $0\leq s\leq t$.
As before, $Y(t)$ is the principal fundamental matrix of $\dot{y}(t) = A(t)y$
and $A(t)$ is continuous.
\label{thm-4-2}
\end{thm}

\begin{proof}
As usual, let the discontinuity of $\yt{t}$ at $t=kh$ be $h\alpha(h)v_k$.
Let $n = \floor{t/h}$ and $h_r=t-nh$.
Then 
$$\yt{kh} = Y(kh)\inv{Y}((k-1)h;h)\yt{(k-1)h} + h\alpha(h)v_k,$$
for $1\leq k\leq n$, and $\yt{t} = Y(t)\inv{Y}(nh) \yt{nh}$. Combine
these equalities to get, 
\begin{equation}
\tilde{y}(t;h) = \sum_{k=1}^{n} Y(t)\inv{Y}(kh) h \alpha(h) v_k,
\label{eqn-s4-1}
\end{equation}
This expression is also used in the proof of Theorem \ref{thm-5-1}.

To prove,
\begin{equation}
\sup_{v(s)} \Norm{\int_0^t Y(t)\inv{Y}(s) v(s) ds} \leq E(t),
\label{eqn-s4-2}
\end{equation}
it is enough to find a $\yt{t}$ for every $v(s)$ such that
\begin{equation}
\int_{0}^{t} Y(t)\inv{Y}(s) v(s) ds = \frac{\yt{t}}{\alpha(h)} + \eta(h),
\label{eqn-s4-3}
\end{equation}
where $\eta(h)\rightarrow 0$ as $h\rightarrow 0$. Take the discontinuity of
$\yt{t}$ at $t=kh$ to be $h\alpha(h)v(kh)$. Use \eqref{eqn-s4-1} and
it is apparent that $\yt{t}/\alpha(h)$ approximates the Riemann integral
$\int_0^t Y(t)\inv{Y}(s) ds$. The standard convergence theorem for
Riemann integrals gives $\eta(h)\rightarrow 0$ as $h\rightarrow 0$.

To show \eqref{eqn-s4-2} in the reverse direction, it is enough to
find a $v(s)$ for a given $\yt{t}$ so that $\eqref{eqn-s4-3}$
holds. First, consider a discontinuous $\tilde{v}(s)$ defined by
$\tilde{v}(s) = v_k$ for $kh \leq s < (k+1)h$.  Clearly,
$\norm{\tilde{v}(s)} \leq 1$.  Then,
\begin{equation}
 \int_0^t Y(t)\inv{Y}(s) \tilde{v}(s) ds = \frac{\yt{t}}{\alpha(h)} +
\delta(h) t,
\label{eqn-s4-4}
\end{equation}
where $\delta(h) = \sup_{0\leq s \leq t} \norm{Y(t)\inv{Y}(s+h) -
Y(t)\inv{Y}(s)}$. Since $Y(s)$ is continuous, $\delta(h)\rightarrow 0$
as $h\rightarrow 0$.
Lusin's theorem \cite{Rudin} is a basic result for approximating measurable
functions by continuous functions. It guarantees a continuous $v(s)$ such
that $\norm{v(s)} \leq \norm{\tilde{v}(s)}$ and $\mu(v(s)\neq \tilde{v}(s)) < h$,
where $\mu$ is the Lebesgue measure. 
If $M = \sup_{0\leq s\leq t}
\norm{Y(t)\inv{Y}(s)}$, then
\begin{equation}
\int_0^t Y(t)\inv{Y}(s) v(s) ds = \int_0^t Y(t)\inv{Y}(s) \tilde{v}(s) ds
+ 2mh,\quad \abs{m} \leq 2M.
\label{eqn-s4-5}
\end{equation}
Together, \eqref{eqn-s4-4} and \eqref{eqn-s4-5} complete the proof.
\end{proof}

\begin{cor}
$$E(t) \leq \int_0^t \norm{Y(t)\inv{Y}(s)} ds.$$
\label{cor-4-1}
\end{cor}

\begin{cor}
For $0 \leq \delta \leq t$ , 
$$E(t) \geq \norm{\int_0^\delta Y(t)\inv{Y}(s)ds}.$$
\label{cor-4-2}
\end{cor}

The definitions of uniform stability and uniform asymptotic stability that
follow seem to have been introduced by Malkin \cite{Malkin}. Theorems
which deduce the stability of a nonlinear system from its linear first
approximation usually (always?) assume the linear first approximation
to be uniformly stable or uniformly asymptotically stable \cite{SC}. The
uniformity assumptions are not explicitly stated sometimes, for example
in \cite{Bellman} and \cite{CL}. In these cases, the $A(t)$ in 
$\dot{y}(t) = A(t)y$ is either constant or periodic, which means that
stability implies uniform stability and asymptotic stability implies
uniform asymptotic stability. Uniformity assumptions are natural in
the theory of Lyapunov functions as well \cite{Yoshizawa}. We will find
them useful for bounding $E(t)$.

\begin{defn}
The solution $x(t;x_0)$ of $\dot{x}(t) = f(t,x)$ is {\it uniformly stable}
if for every $\epsilon > 0$ there exists a $\delta > 0$ such that
$\norm{x(\tau;x_0) - x_\tau'}< \delta$ for $\tau\geq 0$ implies 
$\norm{x(t;x_0) - x(t; \tau, x_\tau')} < \epsilon$ for $t\geq \tau$.
\end{defn}

\begin{defn}
The solution $x(t;x_0)$ is {\it uniformly asymptotically stable} if 
it is uniformly stable and the choice of $\delta$ in the previous definition
can be made in such a way that $\norm{x(t;x_0)-x(t;\tau,x_\tau')}
\rightarrow 0$ as $\tau \rightarrow \infty$ in a uniform way; \ie\,
given $\epsilon'>0$ there exists $T_{\epsilon'}$ such that 
$\norm{x(t;x_0) - x(t;\tau,x_\tau')}<\epsilon'$ for all 
$t > \tau+T_{\epsilon'}$
and $x_\tau'$ satisfying $\norm{x(\tau;x_0) - x_\tau'} < \delta$.
\end{defn}

\begin{thm}
If the zero solution of $\dot{y}(t) = A(t)y$, $y(0) = 0$, is uniformly 
stable, its $E(t)$ is linearly bounded; \ie, $E(t) \leq Kt$ for some
$K>0$ and $0\leq t < \infty$.
\label{thm-4-3}
\end{thm}
\begin{proof}
Uniform stability of the zero solution is equivalent to boundedness of
$\norm{Y(t)\inv{Y}(s)}$ for $t\geq s\geq 0$ \cite{SC}
\cite{Yoshizawa}. If $\norm{Y(t)\inv{Y}(s)}\leq K$ for $t\geq s\geq 0$,
Corollary \ref{cor-4-1} implies $E(t) \leq Kt$.
\end{proof}

\begin{thm}
If the zero solution of $\dot{y}(t) = A(t)y$, $y(0)=0$, is uniformly
asymptotically stable, its $E(t)$ is bounded by a constant; \ie\, $E(t)<K$
for some $K>0$ and $0\leq t<\infty$.
\label{thm-4-4}
\end{thm}
\begin{proof}
Uniform asymptotic stability of the zero solution is equivalent to
$\norm{Y(t)\inv{Y}(s)} < Me^{-\nu (t-s)}$ for $\nu > 0$, $M > 0$, and
$t\geq s \geq 0$ \cite{SC} \cite{Yoshizawa}.
Again, we can use Corollary \ref{cor-4-1} to complete the proof.
\end{proof}

Theorem \ref{thm-4-3} implies that $E(t)$ for the solution of $\dot{x}(t)
=Ax$, $x(0) = x_0$, is linearly bounded if all the eigenvalues of $A$
have negative or zero real parts, and the ones with zero real part are
simple. If all the eigenvalues of $A$ have strictly negative real parts
then, in fact, $E(t)$ is bounded by a constant by Theorem \ref{thm-4-4}.
The necessary stability properties of the zero solution of $\dot{x}(t)=Ax$
are verified in numerous places including \cite{SC}.

%% file: s5.ltx
\section{$E(t)$ for Nonlinear Systems}

This section gives two approaches to the analysis of $E(t)$ of nonlinear
systems. Theorem \ref{thm-5-1} proves that the $E(t)$ of the solution
of a nonlinear system and the $E(t)$ of the zero solution of its
first approximation are the same. So one approach is to look at the
linearized problem. The other approach is to directly make stability
assumptions about the solution of the nonlinear system (Theorems \ref{thm-5-3}
and \ref{thm-5-4}). Both these approaches are illustrated in Section 7.

\begin{lem}
Assume that $F(t,0)\equiv 0$, that $F(t,y)$ has continuous first order
partial derivatives with respect to $y$, and that $F(t,y)$ is 
continuous with respect to $t$, $t\geq 0$. Then,
\begin{equation*}
F(t,y) = A(t)y + g(t,y),
\end{equation*}
where $A(t) = \frac{\partial F(t,y)}{\partial y}\bigl\lvert_{y=0}$ and
$g(t,y) = o(\norm{y})$ as $y\rightarrow 0$ uniformly over compact
intervals of $t$.
\label{lem-5-1}
\end{lem}
\begin{proof}
See \cite{CL}.
\end{proof}

As  noted in the introduction, the following theorem is implicit in the
work of Henrici \cite{Henrici1} \cite{Henrici2}. It might be possible to
generalize it to partial differential equations. 

\begin{thm}
As in Lemma \ref{lem-5-1}, let $F(t,0) \equiv 0$, let $F(t,y)$ have
continuous first order partial derivatives in $y$ and  
be continuous in $t$. Then the zero solution of 
$$\dot{y}(t) = F(t,y),\quad y(0)=0,$$ 
and the zero solution of
$$\dot{y}(t) = A(t)y,\quad y(0)=0,$$
where $A(t) = \frac{\partial F(t,y)}{\partial y}\bigl\lvert_{y=0}$ have
the same $E(t)$.
\label{thm-5-1}
\end{thm}
\begin{proof}
Let $F(t,y)= A(t)y+g(t,y)$ as in Lemma \ref{lem-5-1}. If needed, take
$h < h_0(t,r)$ for some $r>0$. To reduce clutter in the proof, denote
$\yt{kh}$ by $\tilde{y}_k$ and $Y(kh)$ by $Y_k$; $Y(t)$ is the principal
fundamental matrix of $\dot{y}(t) = A(t)y$. Let $n=\floor{t/h}$
and $h_r = t-nh$.

Routine application of the variation of constants formula \cite{HNW}
\cite{SC} gives
\begin{equation*}
\tilde{y}_k = Y_k \inv{Y}_{k-1} \tilde{y}_{k-1} + \int_{(k-1)h}^{kh}
Y_k \inv{Y}(s)g(s,\tilde{y}(s;h))\, ds + h\alpha(h)v_k,
\end{equation*}
for $k = 1,2,\ldots,n$, and 
\begin{equation*}
\yt{t} = Y(t)\inv{Y}_n \tilde{y}_n+\int_{nh}^t Y(t)\inv{Y}(s)g(s,\yt{s})\, ds.
\end{equation*}
Here, $\tilde{y}_0 = 0$ and $Y_k\inv{Y}(s) = I + O(h)$ for $(k-1)h\leq
s\leq kh$.
Further, $\norm{\tilde{y}(s;h)}\leq te^{Lt} \alpha(h)$ for
$0\leq s\leq t$ by Lemma \ref{lem-2-1}, and $\norm{g(s,y)} = o(y)$ as
$y\rightarrow 0$ uniformly over $s\in[0,t]$. Therefore, we can write
$Y_k\inv{Y}(s)g(s,\yt{s}) = C\alpha(h)\eta(h)u(s)$ 
for $0\leq s \leq t$, a constant
$C$ which can depend on $t$ but not on $h$, and $\norm{u(s)} \leq 1$,
where $u(s)\in R^d$; here $\eta(h)\rightarrow 0$ as $h\rightarrow 0$.
The expressions for $\tilde{y}_k$ and $\yt{t}$ become,
\begin{align*}
\tilde{y}_k &= Y_k \inv{Y}_{k-1} \tilde{y}_{k-1} 
+ C\alpha(h)\eta(h) \int_{(k-1)h}^{kh} u(s) ds + h \alpha(h) v_k, \\
\yt{t} &= Y(t) \inv{Y}_n \tilde{y}_n + C\alpha(h)\eta(h) \int_{nh}^t u(s)ds.
\end{align*}
Combining these equalities,
\begin{equation*}
\begin{split}
\frac{\yt{t}}{\alpha(h)} = h \sum_{k=1}^n Y(t) \inv{Y}_k v_k 
&+ C h \eta(h) \sum_{k=1}^n Y(t) \inv{Y}_k \int_{(k-1)h}^{kh}u(s) ds \\
&+ C\eta(h) \int_{nh}^t u(s) ds.
\end{split}
\end{equation*}
Both the $2$nd and $3$rd terms vanish as $h\rightarrow 0$ because
of $\eta(h)$. Thus the nonlinear term $g(t,y)$ has no effect on
$E(t)$. The proof can be formally completed using \eqref{eqn-s4-1}.
\end{proof}

\begin{cor}
Let $f(t,x)$ be continuous in $t$ and have continuous first order partial
derivatives with respect to $x$. The $E(t)$ of the solution $x(t;x_0)$ of
$\dot{x}(t) = f(t,x)$, $x(0) = x_0$ , and  the $E(t)$ of the zero
solution of $\dot{y}(t) = A(t)y$, $y(0)=0$, where 
$A(t) = \frac{\partial f(t,x)}{\partial x}\bigl\lvert_{x = x(t;x_0)}$,
are the same. In fact,
\begin{equation*}
E(t) = \sup_{v(s)} \Norm{\int_0^t \frac{\partial x(t)}{\partial x(s)} v(s) ds},
\end{equation*}
where the supremum is taken over continuous functions $v(s)$ satisfying
$\norm{v(s)} \leq 1$.
\label{cor-5-0}
\end{cor}
\begin{proof}
For a proof of the above formula for $E(t)$, note that if $Y(t)$ is
the principal fundamental matrix of $\dot{y}(t)=A(t)y$, then
$Y(t)\inv{Y}(s) = \frac{\partial x(t)}{\partial x(s)}$; for a proof
see any basic book on ODEs. Plug this into Theorem \ref{thm-4-2} to
get the formula. 
\end{proof}

For comparison with Theorem \ref{thm-5-1}, we state Theorem 33 of Chapter IX
of Sansone and Conti \cite{SC}. This result is typical of Lyapunov's theory
of first approximation.

\begin{thm}
Assume the zero solution of $\dot{y}(t) = A(t) y$, $t\geq 0$,
is uniformly asymptotically stable. If $g(t,y)$ is continuous in $t$ and
$y$, and $g(t,y) = o(\norm{y})$ uniformly over $t\geq 0$ as $y\rightarrow 0$,
then the zero solution of $\dot{y}(t) = A(t)y + g(t,y)$ is also 
uniformly asymptotically stable.
\label{thm-5-2}
\end{thm}

The assumption about $g(t,y)$ being $o(\norm{y})$ uniformly over the
semi-infinite interval $t\geq 0$ in Theorem \ref{thm-5-2} is quite
stringent. A nonlinearity of the form $g(t,y) = t[y_1^2,\ldots,y_d^2]^T$
does not satisfy that assumption. But the theorem does not hold if the
assumption is weakened to what is known about $g(t,y)$ from 
Lemma \ref{lem-5-1}; for a counterexample, see Bellman \cite[p. 87]{Bellman}.
A stringent assumption about $g(t,y)$ is not required in Theorem \ref{thm-5-1}
because how small an $h$ we take to get convergence of $e(t;h)/\alpha(h)$
to $E(t)$ can depend upon $t$. In definitions of stability, in contrast,
we need one small perturbation which stays ``close'' till $t=\infty$.

The rest of this section is about bounding $E(t)$ by making stability 
assumptions about the solution of a nonlinear system. For linear systems,
this program was carried out in Section 4. We introduce a technique for
error analysis of approximate solutions which uses Lyapunov functions.
But first an example to point out some difficulties.

\begin{exmp}
Let us first consider the zero solution of $\dot{y}(t) = y - e^{t}y^{3}$,
$y(0) = 0$, $t\geq 0$. Its $E(t)$, by Theorems \ref{thm-5-1} and 
\ref{thm-4-1}, is $e^t-1$. But we show that the zero solution is actually
uniformly asymptotically stable. It is even exponentially asymptotically
stable.

\begin{figure}
\centerline{
\psfig{figure=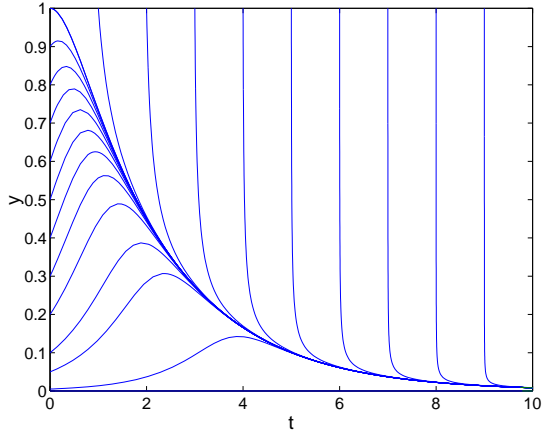,height=2.3in}}
\caption[xyz]{ The portrait of trajectories of $\dot{y}(t) = y -
e^ty^3$. All solutions tend towards the curve $\dot{y}(t) =
e^{-t/2}$.}
\label{fig-3}
\end{figure}

Figure \ref{fig-3} is the portrait of trajectories of $\dot{y}(\tau) = y - e^ty^3$. The
portrait for $y\leq 0$ is a reflection about $y=0$. Thus we 
can restrict ourselves to trajectories which are always in the upper half
plane. Every point $(t,y)$ with $y\geq e^{-t/2}$ is on a trajectory
that is pointed downwards. 

We now verify uniform stability of the zero solution using the last
observation. Given $\epsilon > 0$, choose $t_\epsilon$ so that
$e^{-t_\epsilon/2} < \epsilon$. By continuity properties, we can
choose a $\delta > 0$ so that if $\tau \leq t_\epsilon$ and
$y_\tau\leq \delta$, the trajectory through $(\tau, y_\tau)$ stays
below $\epsilon$ till $t_\epsilon$. The maximum possible height (along $y$)
of a trajectory beginning at $(\tau, y_\tau)$, $\tau\geq t_\epsilon$
and $y_\tau < \delta$, is bounded by the larger of $e^{-\tau/2}$ and
$\abs{y_\tau}$. Since $e^{-\tau/2} < \epsilon $ and $\delta < \epsilon$,
uniform stability is verified.

The verification of uniform asymptotic stability will be sketchy. We base
it on the following facts:
\begin{enumerate}
\item The solution of $\dot{y}(t) = y-e^ty^3$, $y(0) = y_0$ tends to zero
as $t\rightarrow \infty$ for $0\leq y_0\leq 1$,
\item Further, $y(t;y_0)\leq y(t;1)$ for $t\geq 0$, $0\leq y_0\leq 1$,
\item And, $0\leq y(t+\tau; \tau; y_0)\leq y(t;0;y_0)$ for any $y_0\geq 0$,
$\tau \geq 0$, $t\geq 0$.
\end{enumerate}
The proof of item 1 involves a bit of elementary work which we do at
the end. Item 2 is trivial. For item 3, think of $y(t;0,y_0)$ as the
solution of $\dot{y}(t)=y-e^ty^3$, $y(0)=y_0$, and of $y(t+\tau; \tau,
y_0) = z(t)$ as the solution of $\dot{z}(t) = z - e^{t+\tau} z^3$,
$z(0) = y_0$, and use a differential inequality \cite[p. 27]{Hartman}.

Now let $T_\epsilon$ be such that $\norm{y(t;y_0)} <\epsilon$ for $t\geq
T_\epsilon$ and $y_0=1$. Then $\norm{y(t+\tau; \tau, y_0)} < \epsilon$
for any $\tau\geq 0$, $t \geq T_\epsilon$ and $\abs{y_0} \leq 1$. Thus,
the zero solution is uniformly asymptotically stable.

To prove item 1, it is enough to verify that the solution with the
initial condition $y(0) = 1$ satisfies $y(t) < 2e^{-t/2}$ for $t\geq 0$.
This is obvious from the portrait of trajectories: every trajectory that
cuts the curve $y = 2e^{-t/2}$ is in the downward direction; therefore
any trajectory that starts below that curve has to stay below it for
$t\geq 0$.  In fact, $y(t) < 2e^{-t/2}$ for the trajectory with $y(0)=1$
implies that the zero solution is {\it exponentially asymptotically
stable}, if we adopt Yoshizawa's definition of exponential
asymptotic stability \cite{Yoshizawa}.

\label{exmp-5-1}
\end{exmp}

Here we have an example which is uniformly asymptotically stable, yet has
an $E(t)$ which increases exponentially with $t$. This is
possible because
the approximations $\yt{t}$ can introduce errors at all points $kh$ 
while the definitions of stability allow just one perturbation.

The next two theorems are nonlinear analogues of Theorem \ref{thm-4-3}
and \ref{thm-4-4}. The proofs this time rely on the theory of Lyapunov
functions. Let us note that Lyapunov functions $V(t,y)$ are always assumed
to be continuous in $t$ and locally Lipshitz in $x$.
$V_F'(t,y)$ is the rate of increase of $V(t,y)$ along a solution
of $\dot{y}(t) = F(t,y)$ which goes through $(t,y)$. More precisely,
if $y(t+\tau;t,y)$ is such a solution, then
\begin{equation*}
V_F'(t,y) = \limsup_{\delta\rightarrow 0^{+}}
\frac{V(t+\delta, y(t+\delta; t,y)) - V(t,y)}{\delta}.
\end{equation*}
If $V_F'(t,y)\leq a V(t,y)$ then $V(t+\delta, y(t+\delta; y,t)) \leq 
e^{a\delta} V(t,y)$, and if $V_F'(t,y) \leq 0$ then 
$V(t+\delta, y(t+\delta; y,t)) \leq  V(t,y)$; these two facts can
be inferred from differential inequalities \cite{HNW} \cite{Hartman}.

We say that $F(t,y)$ is uniformly Lipshitz in a neighbourhood of zero
if $\norm{F(t,y_1) - F(t,y_2)} < L\norm{y_1- y_2}$, for a constant
$L>0$ and any $y_i$ with $\norm{y_1} < r$, $\norm{y_2} < r$ where
$r>0$; $L$ is the same constant for any $t\geq 0$.

\begin{thm}
Let $F(t,y)$ be uniformly Lipshitz in $y$ in a neighbourhood of
$0$. Assume that the zero solution of $\dot{y}(t) = F(t,y)$, $y(0)=0$,
is exponentially stable in the sense that
\begin{equation*}
\norm{y(t; t_0, y_0)} \leq K e^{-c(t-t_0)} \norm{y_0}
\end{equation*}
for $\norm{y_0} < r$, $t_0\geq 0$, $t\geq t_0$, $c>0$, and $K>0$. Then
$E(t)$ of the zero solution is bounded above by a constant.
\label{thm-5-3}
\end{thm}
\begin{proof}
Stability assumptions in the theorem imply existence of Lyapunov
function with following properties (Yoshizawa \cite[p. 97]{Yoshizawa},
Hale \cite{Hale}):
\begin{enumerate}
\item $\norm{y} \leq V(t,y) \leq C \norm{y}$, where $C>0$ is a constant,
\item $\abs{V(t,y_1) - V(t,y_2)} \leq L \norm{y_1-y_2}$,
\item $V_F'(t,y) \leq -q c V(t,x)$ for some $0<q<1$.
\end{enumerate}
The domain of $V(t,y)$ is $t\geq 0$ and $\norm{y}\leq r$ for $r>0$.

Take $h < h_0(t,r)$. Because of item 3, $V(t,\yt{t})$ decreases at
least by a factor $e^{-qch}$ along the approximate solution when $t$
increases from $kh$ to $(k+1)h$; on that interval of $t$ the
approximate solution follows the exact solution till the discontinuity
at $t=(k+1)h$. The discontinuity can cause an increase in $V(t,y)$ of
at most $L$ times its magnitude by item 2. Therefore,
\begin{equation*}
V(kh; \yt{kh}) \leq e^{-qch}V((k-1)h, \yt{(k-1)h}) + Lh\alpha(h),
\end{equation*}
for $k=1,2,\ldots,n$ and
\begin{equation*}
V(t,\yt{t}) \leq e^{-qch_r} V(nh, \yt{nh}).
\end{equation*}
Combining these inequalities, we get 
\begin{align*}
V(t,\yt{t}) &\leq e^{-qch_r}L\alpha(h) \Bigl(\frac{1-e^{-nqch}}{1-e^{-qch}}
\Bigr) \\
&\leq K \alpha(h).
\end{align*}
That $K$ above can be a constant independent of $h$ and $n$ can be
deduced from basic calculus using $h<h_0$. Now, by item 1, 
$\norm{\yt{t}} \leq C \alpha(h)$ implying a constant upper bound
for $E(t)$.
\end{proof}

The difficulty in proving Theorem \ref{thm-5-3} directly using the
norm $\norm{\cdot}$ is that when $K>1$ the  discretization error
might actually be amplified by a factor greater than $1$ over any
given time step. Since we have to make the worst possible assumption
at every time step, the final bound on $E(t)$ obtained this way
will actually be exponential in $t$ when $K > 1$. The proof of
Theorem \ref{thm-5-3} uses a carefully constructed Lyapunov function
to get around this difficulty.  

\begin{thm}
Assume as in the previous theorem that $F(t,y)$ is uniformly Lipshitz 
with respect to $t$ in a neighbourhood of $y=0$. If the zero solution
of $\dot{y}(t) = F(t,y)$, $y(0)=0$, $t\geq 0$, is uniformly asymptotically
stable, then $E(t) \leq Kt$ for some constant $K$.
\label{thm-5-4}
\end{thm}
\begin{proof}
Stability assumptions in this theorem imply the existence of a Lyapunov
function with the following properties:
(Hale \cite[Theorem 4.2, Chapter X]{Hale}, Yoshizawa \cite{Yoshizawa})
\begin{enumerate}
\item $\norm{y} \leq V(t,y)$,
\item $V(t,0) \equiv 0$,
\item $V_F'(t,y) \leq 0$, where $V_F'(t,y)$ is defined as in the previous
proof,
\item $\abs{V(t,y_1) - V(t,y_2)} \leq K \norm{y_1-y_2}$ for some constant
$K>0$.
\end{enumerate}
The domain of definition of $V(t,y)$ is the same as in the previous proof.

The proof is similar to that of Theorem \ref{thm-5-3}, but this time
\begin{equation*}
V(kh, \yt{kh}) \leq V((k-1)h, \yt{(k-1)h}) + Kh\alpha(h)
\end{equation*}
for $k=0,1,\ldots,n-1$
and 
$$V(t,\yt{t}) \leq V(nh, \yt{nh}).$$
Combining these inequalities, we
have $V(t,\yt{t})\leq Kt\alpha(h)$.  As before, $\norm{\yt{t}}\leq
Kt\alpha(h)$, which this time implies that $E(t)\leq Kt$.
\end{proof}

\begin{cor}
Let $x(t;x_0)$ be the nonzero solution of the initial value problem
$\dot{x}(t) = f(t,x)$, $x(0)=x_0$. Assume $f(t,x)$ is uniformly Lipshitz
in $x$ in a neighbourhood of $x(t;x_0)$. Then uniform asymptotic
stability of $x(t;x_0)$ implies that its $E(t)$ is linearly bounded.
\label{cor-5-1}
\end{cor}

Let us call attention to the necessity of making a uniform Lipshitz 
assumption about $F(t,x)$ or $f(t,x)$ in Theorem \ref{thm-5-4} or Corollary
\ref{cor-5-1}. Example \ref{exmp-5-1}, which is uniformly asymptotically
stable, does not satisfy the uniform Lipshitz assumption and has an
$E(t)$ which increases exponentially. We do not know if Theorem \ref{thm-5-4}
is still true if the assumption of uniform asymptotic stability is 
weakened to just uniform stability. If such a theorem were true, its
wider applicability might be of use. Table \ref{table-1} summarizes all
of Sections 5 and 6 except Theorem \ref{thm-5-1}. 
\begin{table}
\begin{center}
\begin{tabular}{|p{1.30in}|p{2.25in}|p{2.25in}|}
\hline & Stability & $E(t)$ can increase exponentially even if $\norm{A(t)}$
is bounded \\
\cline{2-3} & Asymptotic stability & $E(t)$ can increase exponentially even
if $\norm{A(t)}$ is bounded; Example \ref{exmp-4-3} \\
\cline{2-3} Linear Problems & Uniform stability & $E(t)$ must be linearly
bounded; Theorem \ref{thm-4-3}\\
\cline{2-3} & Uniform asymptotic stability & $E(t)$ must be bounded by
a constant; Theorem \ref{thm-4-4}\\
\hline & Uniform stability & $E(t)$ can increase exponentially \\
\cline{2-3} &Uniform stability with  uniform Lipshitz assumption &
           Not known if $E(t)$ must be linearly bounded \\
\cline{2-3} &Uniform asymptotic stability & $E(t)$ can increase 
exponentially; Example \ref{exmp-5-1} \\
\cline{2-3} Nonlinear problems & Uniform asymptotic stability with
uniform Lipshitz assumption & $E(t)$ must be linearly bounded;
Theorem \ref{thm-5-4} \\
\cline{2-3} & Exponential stability as in Theorem \ref{thm-5-3} &
Not known if $E(t)$ must be linearly bounded \\
\cline{2-3} & Exponential stability as in Theorem \ref{thm-5-3} with
uniform Lipshitz assumption & $E(t)$ must be bounded by a constant;
Theorem \ref{thm-5-3}\\
\hline
\end{tabular}
\end{center}
\caption[xyz]{Summary of part of Sections 5 and 6. The second column
is the stability assumption about the solution; the last column says
what is known about the conditioning function $E(t)$ corresponding to
such a solution.}
\label{table-1}
\end{table}

%% file: s6.ltx
\section{Three Applications to Dynamical Systems}
We give three examples to illustrate the applicability of our methods
for bounding the accumulation of global error.  

\subsection{Hyperbolic Sinks of $C^{1}$ Dynamical Systems}

Let $p$ be a fixed point of a $C^1$ dynamical system $\dot{x}(t) =
f(x)$; \ie\, let $f(p)=0$. Then $p$ is a hyperbolic sink, if all the
eigenvalues of $\frac{\partial f}{\partial x}\bigl\lvert_{x=p}$ have
strictly negative real parts. The following theorem can be derived
using Chapter 6 of \cite{SH}. We give the theorem here because our
method of proof is different.

\begin{thm}
Let $x(t;x_0)$ be a trajectory of the dynamical system $\dot{x}(t) =
f(x)$, $f\in C^1(R^d)$, which falls into a hyperbolic sink $p$ as
$t\rightarrow \infty$. Then its $E(t)$ is bounded above by a constant.
\label{thm-6-1}
\end{thm}

\begin{proof}
Without loss of generality, take $p=0$. By \cite[p. 150]{Robinson}, there
is a neighbourhood $U_0$ of $0$ such that $x_0 \in U_0$ implies
\begin{equation*}
\norm{x(t;x_0)} < c e^{-at}
\end{equation*}
for constants $a>0$, $c>0$, and for $t\geq 0$.

Using continuity properties of solutions of differential equations
(or openness of the basin of attraction), we can assume an open 
neighbourhood $U$ of $\{x(t;x_0) \bigl\lvert t\geq 0\}$ in $R^d$ and
a compact set $K$ containing $U$ such that $K$ and consequently $U$ are
both contained in the basin of attraction of $p=0$. Since all trajectories
beginning in $K$ enter $U_0$ in a finite amount of time, which depends
only on $K$, we can assume
\begin{equation*}
\norm{x(t;x_0)} < C e^{-at}
\end{equation*}
for constants $a>0$, $C>0$, for $t\geq 0$, and for any $x_0\in K$.
We can also take $K$ and $U$ to be invariant under the flow.

For any trajectory $x(t;x_0)$, $t\geq 0$, with $x_0\in U$, there 
exists $r>0$ such that if $\norm{x(t;x_0)-x_1} < r$, $t\geq 0$ then
$x_1\in K$. 
For such an $x_0$, if $\norm{\delta} < r$,
\begin{equation*}
\norm{x(t;x_0) - x(t; \tau, x(\tau;x_0) + \delta)} < 2Ce^{-at}
\end{equation*}
for $t\geq \tau\geq 0$. If $x_0\in\!\!\!\!\!\!/ U$ but is in the basin of
attraction of $p$, the trajectory $x(t;x_0)$ enters $U$ in a finite
amount of time. So we can get the same kind of bound as above by
adjusting $C$ and $a$ if necessary.

To apply Theorem \ref{thm-4-3} and deduce boundednes of $E(t)$, we need
only verify uniformity of the Lipshitz condition on $f(x)$ for 
$x\in K$. This is trivial since $K$ is compact and $f(x)$ is $C^1$.

\end{proof}

\subsection{Hyperbolic, Attracting Cycles of $C^{1+\epsilon}$ Dynamical Systems}

Let $x(t)$, $t\geq 0$, be a periodic orbit of the $C^1$ dynamical
system $\dot{x}(t) = f(x)$ in $R^d$. Let $T>0$ be its period so that
$x(t+T) = x(t)$.  Denote the set of points on this orbit by $\gamma$.

The characteristic multipliers of the cycle $\gamma$ can be defined in
two ways. One is to pick a point $p\in\gamma$, take a cross-section $\Sigma$
at $p$, define a Poincar\'{e} map for $\Sigma$, and then define the
characteristic multipliers as the $(d-1)$ eigenvalues of the linearization
of the Poincar\'{e} map at $p$. The other way is to consider the linear
first approximation $\dot{y}(t) = A(t)y$ on the cycle $\gamma$. Obviously,
$A(t+T) = A(t)$ for $t\geq 0$. The Floquet numbers of this linear system 
can also be used to define characteristic multipliers. For a lucid
account of these matters, see Robinson \cite{Robinson}.

The cycle $\gamma$ is hyperbolic and attracting if all its characteristic
multipliers are strictly less than $1$ in magnitude.

\begin{thm}
Let $x(t;x_0)$, $t\geq 0$, be an orbit of a $C^{1+\epsilon}$ dynamical system
$\dot{x}(t) = f(x)$ in $R^d$ which falls into a hyperbolic, attracting
cycle $\gamma$ as $t\rightarrow \infty$. Then its $E(t)$ is linearly
bounded from above.
\label{thm-6-2}
\end{thm}

Let us first prove the following lemma. If any one solution of a linear
system is uniformly stable, so is every other solution. So we might
speak of the linear system itself as being uniformly stable.

\begin{lem}
Assume $x_0\in \gamma$ so that $x(t;x_0)$ is a periodic orbit. Let its
linear first approximation be $\dot{y}(t) = A(t)y$, $t\geq 0$. If $\gamma$
is hyperbolic and attracting,  $\dot{y}(t) = A(t)y$ is uniformly
stable, and the $E(t)$ associated with $x(t;x_0)$ is linearly bounded.
\label{lem-6-2-1}
\end{lem}
\begin{proof}
Uniform stability of $\dot{y}(t) = A(t)y$ is an easy consequence of the
characteristic multipliers of $\gamma$ being strictly less than $1$.
See Chapter IX of \cite{SC}. The linear bound on $E(t)$ is implied by
Theorem \ref{thm-4-3} and Corollary \ref{cor-5-1}. 
\end{proof}

Lemma \ref{lem-6-2-1} is contained in a different form in the work of
Cano and Sanz-Serna \cite{CanoS}. But let us note that Theorem \ref{thm-6-2}
goes beyond Lemma \ref{lem-6-2-1} in a significant way. In practice,
it is highly unlikely that $x_0$ itself is on the cycle $\gamma$. But it
is often easy to find $x_0$ so that $x(t;x_0)$ falls into a cycle 
$\gamma$.

The following lemma is known as the Dini-Hukuhara-Caligo theorem. It
is Corollary 1 of Chapter IX of \cite{SC}. Its proof, which we omit,
is short and simple, and illustrative of an important technique in
stability theory.

\begin{lem}
Assume the linear system $\dot{y}(t) = A(t)y$, $t\geq 0$, is uniformly
stable. Assume also that $B(t)$, $t\geq 0$, is continuous with
$\int_{0}^{\infty} \norm{B(t)} dt < \infty$. Then the linear system
$\dot{y}(t) = (A(t)+B(t)) y$ is also uniformly stable.
\label{lem-6-2-2}
\end{lem}

\begin{proof}[Proof of Theorem \ref{thm-6-2}]
By Hartman \cite[p. 254]{Hartman}, there exists a point $x_0'\in \gamma$
such that
\begin{equation}
\norm{x(t;x_0) - x(t;x_0')} < c e^{-at},
\label{eqn-6-local-1}
\end{equation}
for constants $a>0$, $c>0$, and for $t\geq 0$. This is called convergence
in phase \cite{Robinson}.

Let $\dot{y}(t) = A(t)y$, where
$A(t) = \frac{\partial f}{\partial x}\bigl\lvert_{x=x(t;x_0')}$, be the
first approximation along $x(t;x_0)$. By Lemma \ref{lem-6-2-1}, this
linear system is uniformly stable.

Let $\dot{y}(t) = (A(t) + B(t))y$, where $A(t)+B(t) = \frac{\partial
f}{\partial x}\bigl\lvert_{x=x(t;x_0)}$ be the first approximation
along $x(t;x_0)$. The estimate \eqref{eqn-6-local-1} for convergence
in phase implies
\begin{equation*}
\norm{B(t)} < c_1 e^{-a_1t},
\end{equation*}
for constant $a_1>0$ and $c_1>0$. This is because both $x(t;x_0)$ and
$x(t;x_0')$ stay within a compact region of $R^d$, and $f(x)$ is
$C^{1+\epsilon}$
(this where we need the $\epsilon$ in $C^{1+\epsilon}$).
By Lemma \ref{lem-6-2-2}, $\dot{y}(t) =
(A(t)+B(t))y$ is also uniformly stable.

Since $A(t)+B(t)$ gives the linearization of $x(t;x_0)$,
the proof is easy to complete using Theorem \ref{thm-4-3} and
Corollary \ref{cor-5-1}.
\end{proof}

\subsection{Normally Contracting Manifolds with Quasiperiodic Flows}

Let us introduce the notation $\phi_t$ for the flow induced on $R^d$
by $\dot{x}(t) = f(x)$.  With this notation $x(t;x_0) = \phi_t x_0$.
Let $V$ be a compact $C^1$ manifold which is invariant under this
flow. We consider the situation when the flow on $V$ is {\it
differentiably conjugate} to quasiperiodic flow on a torus, and $V$ is
normally hyperbolic and contracting, or briefly, {\it normally
contracting}.  We now explain the two italicized concepts in this
paragraph.

A torus $T^n$ is $n$ copies of the circle $S^1$. If the angle on the
$i$th circle is parameterized by $\theta_i$, a quasiperiodic flow 
on $T^n$ is of the form $\theta_i(t) = (\theta_i(0) + \alpha_i t) \mod 2\pi$.
In fact, the flow is periodic if the $\alpha_i$ are all mutually 
commensurable. Denote this flow by $\psi_t$.

When we say that the flow $\phi_t$ is differentiably conjugate to
quasiperiodic flow on a torus, we mean that there exists a $C^1$ 
homeomorphism $h:V\rightarrow T^n$ such that $h(\phi_t x) = \psi_t(h x)$
for $x\in V$ and $t\geq 0$. 

To define normal contractivity \cite{HPS} \cite{Robinson}, associate a
direct sum decomposition $T_x \oplus N_x$ of $R^d$ with every $x$ in
$V$. In this splitting $T_x$ is the tangent space of $V$ at $x$, and
$N_x$, the normal space, varies continuously with $x$. If the $N_x$
can be chosen so that
\begin{equation*}
\norm{\Pi_{N_y} \frac{\partial y}{\partial x} \bigl\rvert N_x} < ce^{-\mu t},
\end{equation*}
where $y = \phi_t x$, the matrix inside the norm is the restriction
of the derivative $\frac{\partial y}{\partial x}$ to
act from $N_x$ to $N_y$, $c>0$, and $\mu >0$, then $V$ is 
normally contracting. Usually, the definition of normal contractivity
comes with an other assumption which says contraction in the normal
direction dominates any contraction on the manifold $V$. But since
we have assumed that the flow on $V$ is differentiably conjugate to
quasiperiodic flow on a torus, this other assumption can be dropped.

Obviously, the tangent spaces $T_x$ are invariant under the derivative
map $\frac{\partial \phi_t x}{\partial x}$. It is actually possible to
choose $N_x$ so that they too are invariant under the derivative map
\cite{HPS} \cite{Robinson}.
We take this to be the case. So the derivative map $\frac{\partial
\phi_t x}{\partial x}$ maps $T_x$ to $T_{\phi_t x}$ and $N_x$ to
$N_{\phi_t x}$.

\begin{thm}
Let $x(t;x_0)$ be a trajectory of the $C^{1+\epsilon}$ flow
$\dot{x}(t) = f(x)$ which falls into a normally contracting and
invariant manifold $V$. Assume that the flow on $V$ is differentiably
conjugate to quasiperiodic flow on a torus. Then the $E(t)$ of
$x(t;x_0)$ is linearly bounded.
\label{thm-6-3}
\end{thm}

The above theorem generalizes Theorem \ref{thm-6-2}. Its proof is exactly
analogous. We begin with a lemma about trajectories that begin on $V$.

\begin{lem}
Let $x_0\in V$ so that $x(t;x_0)$ stays on $V$ for $t\geq 0$.
Make the same assumptions about $V$ as in Theorem \ref{thm-6-3}.
Then the linearization $\dot{y}(t) = A(t)y$ along $x(t;x_0)$ is
uniformly stable, and the $E(t)$ of $x(t;x_0)$ is linearly bounded.
\end{lem}
\begin{proof}
The principal fundamental matrix of the linearization in the lemma is
given by $Y(t) = \frac{\partial \phi_t x_0}{\partial x_0}$, which is
the derivative map. Therefore, it is enough if we show that
$\norm{\frac{\partial \phi_t x}{\partial x}}$ is bounded by a constant
for any $x\in V$ and $t\geq 0$.

We already know that the maps induced by the derivative map between
tangent spaces and between normal spaces are bounded in norm because
of differentiable conjugacy to flow on a torus and normal contractivity,
respectively. Since the tangent spaces and normal spaces are both 
invariant under the derivative map, it is enough if we show that
the angle between $T_x$ and $N_x$ (in the sense of the CS decomposition)
is bounded away from $0$. That this angle is bounded away from $0$ 
is implied by the compactness of $V$.
\end{proof}

The proof of Theorem \ref{thm-6-3} can be completed exactly like the
proof of Theorem \ref{thm-6-2} using the following result about 
convergence in phase.

\begin{thm} 
As in Theorem \ref{thm-6-2}, let $x(t;x_0)$ be a trajectory of the
$C^{1}$ flow $\dot{x}(t) = f(x)$ which falls into a normally
contracting and invariant manifold $V$, and let the flow on $V$ be
differentiably conjugate to quasiperiodic flow on a torus.
Then there exists $x_0'\in V$ such that
\begin{equation*}
\norm{x(t;x_0) - x(t;x_0')} < c e^{-at},
\end{equation*}
for positive constants $c$ and $a$, and $t\geq 0$.
\end{thm}
\begin{proof}
This theorem can be deduced from Theorem 4.1 and the remark following
its proof in Hirsch, Pugh and Shub \cite{HPS}. See in particular
part (a) of that theorem about stable manifolds and part (g) 
about conjugacy to linearized flows.
\end{proof}

%% file: s7.ltx
\section{Three Remarks}
This last section is a collection of parenthetical remarks about
matters which are related to $E(t)$ and which we have not investigated
in detail. 

\begin{description}
\item[(i)] {\it Multistep methods and variable time stepping}. The
model for discretization errors in Section 2 is adapted to single step
methods with constant step sizes. For linear multistep methods with
constant step sizes, we believe the accumulation of global error can
be worse but not better (after excluding some trivial cases) than 
indicated by $E(t)$. 

For variable time stepping, it is sometimes but not always true that
the ratio of the largest to the smallest time step is bounded by
a constant \cite{SN} \cite{Stuart}. In this case at least, the effect
of variable time stepping is to improve global errors by a constant
factor over the indication given by $E(t)$. 

\item[(ii)] {\it One sided Lipshitz conditions}.
Let us briefly summarize
the approach to global error analysis using one sided Lipshitz conditions
for the linear system $\dot{y}(t) = A(t)y$.  Since $\norm{y(t)}^2
= y^T(t) y(t)$, we have
\begin{equation*}
\frac{ d \norm{y(t)}^2}{dt} = y^T(t)(A^T(t) + A(t))y(t).
\end{equation*}
If $\lambda(t)$ is the maximum eigenvalue of $A^T(t) + A(t)$, then
\begin{equation*}
\frac{ d \norm{y(t)}^2}{dt} \leq \lambda(t) \norm{y(t)}^2.
\end{equation*}
Thus upper bounds for $\norm{y(t)}$, and hence for $\norm{Y(t)}$ where
$Y(t)$ is the principal fundamental matrix of the linear system,
can be written down in terms of $\lambda(t)$. These can be plugged
into Theorem \ref{thm-4-2} to get bounds on $E(t)$ and hence the
accumulation of global error. For a detailed account, see \cite{HNW}.

In our view, one sided Lipshitz conditions are basically a way to get
a handle on stability by looking at the evolution of the norm of 
$y(t)$. This is a far less general approach to stability than the
two methods of Lyapunov we have used (this deficiency of norms is
something Lyapunov must have realized). It is also of far lesser
applicability; we do not see a way to derive any of the results in
Section 7 by using one sided Lipshitz conditions.

\item[(iii)] {\it Numerical symplectic integrators and Hamiltonian
problems.} 

Special methods which preserve properties like symplecticity or energy
conservation or volume conservation in phase space show milder
increase of global errors with $t$ than general methods in numerically
integrating Hamiltonian problems. Calvo and Sanz-Serna \cite{CalvoS}
showed using careful numerical experiments that the accumulation
of global error with $t$ was linear for a symplectic method but
quadratic for a Runge-Kutta method when the methods were applied to 
Kepler's problem. A full explanation came after \cite{CH} and \cite{ES}
in the paper \cite{CanoS} by Cano and Sanz-Serna. There are 
numerical experiments in Quispell and Dyt \cite{QD} which are
unexplained. 

The conditioning function $E(t)$ is not suitable for studying the accumulation
of global error in these special methods. The class of local discretization
errors allowed by the model in Section 2 is too general. The discretization
errors made by these special methods are restricted in some way; for
example, the discretization errors of an energy preserving method cannot
take the approximation out of the constant energy manifold. But is it
possible to define  conditioning functions which are adapted to these
special methods? We think it might be possible. The concepts of stability
needed to bound such conditioning fuctions will also, no doubt, be 
different.

\end{description}

\section{Acknowledgements}
It is a pleasure to thank Arieh Iserles and Mike Shub for encouragement
and enlightening discussions. I have benefited from talking to several
visitors to M.S.R.I., Berkeley, in Fall 98 about this work. I thank Timo Eirola
and Andrew Stuart for suggestions which directly influenced the presentation
of this paper.

%% file: 1999-001.bbl
\begin{thebibliography}{44}

\bibitem{Bellman} R. Bellman, {\em Stability Theory of Differential 
Equations}, McGraw-Hill, New York, 1953.

\bibitem{CH} M.P. Calvo and E. Hairer, Accurate long-term integration
of dynamical systems, {\em Appl. Numer. Math. 18} (1995), 95-105.

\bibitem{CalvoS} M.P. Calvo and J.M. Sanz-Serna, The development of
variable-step symplectic integrators, with application to the two-body
problem, {\em SIAM J. Sci. Comput. 14} (1993), 936-952.

\bibitem{CanoS} B. Cano and J.M. Sanz-Serna, Error growth in 
the numerical integration of periodic orbits, with application to
Hamiltonian and reversible systems, {\em SIAM J. Numer. Anal. 34}
(1997), 1391-1417.

\bibitem{CL} E.A. Coddington and N. Levinson, Theory of Ordinary 
Differential Equations, McGraw-Hill, New York, 1955.

\bibitem{Dahlquist} G. Dahlquist, Stability and error bounds in the
numerical integration of ordinary differential equations, 
{\it Trans. of the Royal Inst. of Techn.,} Stockholm, Sweden,
Number 130.

\bibitem{ES} D.J. Estep and A.M. Stuart, The rate of error growth
in Hamiltonian-conserving integrators, {\em Z. Agnew. Math. Phys. 46}
(1995), 407-418.

\bibitem{Gragg} W.B. Gragg, {\em Repeated extrapolation to the limit
in the numerical solution of ordinary differential equations}, Thesis,
Univ. of California, 1964; see also {\em SIAM J. Numer. Anal. 2} (1965),
384-403.  

\bibitem{HNW} E. Hairer, S.P. N{\o}rsett and G. Wanner, {\em Solving
Ordinary Differential Equations I}, Springer-Verlag, New York, 1980.

\bibitem{Hale} J. Hale, {\em Ordinary Differential Equations},
John Wiley, New York, 1969.

\bibitem{Hartman} P. Hartman, {\em Ordinary Differential Equations},
John Wiley, New York, 1973.

\bibitem{Henrici1} P. Henrici, {\em Discrete Variable Methods in Ordinary
Differential Equations}, John Wiley, New York, 1962.

\bibitem{Henrici2} P. Henrici, {\em Error Propogation for Difference
Methods}, John Wiley, New York, 1963.

\bibitem{HPS} M.W. Hirsh, C.C. Pugh and M. Shub, {\it Invariant
Manifolds}, Lecture notes in Mathematics 583, Springer-Verlag,
New York, 1977.
  
\bibitem{IN} A. Iserles and S.P. N{\o}rsett, On the solution of linear
differential equations in Lie groups, {\em Phil. Trans. of Royal Soc. A},
to appear.

\bibitem{IS} A. Iserles and G. S\"{o}derlind, Global bounds on numerical
error for ordinary differential equations, {\it J. Complexity 9} (1993),
97-112.

\bibitem{Lyapunov} A. M. Lyapunov, {\em Probl\`{e}me g\'{e}n\'{e}rale
de la stabilit\'{e} du mouvement},  Comm. Soc. Math. Kharkov 2 (1892),
3 (1893); Ann. Fac. Sci. Toulouse 9 (1907), 204-474; Ann. of Math.
Studies 17, Princeton University Press, Princeton, 1949.

\bibitem{Malkin} J. G. Malkin, {\em Theory of Stability of Motion},
Atomic Energy Commision Tech. Rep. 3352, 1956.

\bibitem{QD} G.R.W. Quispel and C.P. Dyt, Volume-preserving integrators
have linear error growth, {\em Physics Letters A 242} (1998), 25-30.


\bibitem{Robinson} C. Robinson, {\em Dynamical Systems. Stability,
Symbolic Dynamics, and Chaos}, CRC press, Boca Roton, 1995. 

\bibitem{Rudin} W. Rudin, {\em Real and Complex Analysis}, 3rd Edition,
McGraw-Hill, New York, 1987.

\bibitem{SC} G. Sansone and R. Conti, {\em Non-linear Differential
Equations}, Macmillan, New York, 1964.

\bibitem{SN} D. Stoffer and K. Nipp, Invariant curves for variable
step size integrators, {\em BIT 31} (1991), 169-180.

\bibitem{Stuart} A.M. Stuart, Probabilistic and deterministic convergence
proofs for software for initial value problems, 
{\em Numerical Algorithms 14} (1997), 227-260.

\bibitem{SH} A.M. Stuart and A.R. Humphries, {\em Dynamical Systems and
Numerical Analysis}, Cambridge University Press, Cambridge, UK, 1996.

\bibitem{TB} L.N. Trefethen and D. Bau III, {\em Numerical Linear 
Algebra}, SIAM, Philadelphia, 1997.

\bibitem{Yoshizawa} T. Yoshizawa, {\em Stability Theory by Liapunov's 
Second Method}, Mathematical Society of Japan, 1966.

\end{thebibliography}
